\documentclass[10pt,leqno]{article}
\usepackage{amssymb}
\usepackage{hyperref}
\newcommand\qed{\hfill$\sqcap\kern-7.5pt\hbox{$\sqcup$}$}

\newcommand{\Q}{\mathbb{Q}}
\newcommand{\N}{\mathbb{N}}

\newcommand{\R}{\mathbb{R}}

\newtheorem{theo}{Theorem}
\newtheorem{prop}[theo]{Proposition}
\newtheorem{lem}[theo]{Lemma}

\newtheorem{rem}[theo]{Remark}

\newtheorem{defin}[theo]{Definition}

\newcommand{\beqn}{\begin{equation}}
\newcommand{\eeqn}{\end{equation}}
\newcommand{\bear}{\begin{eqnarray}}
\newcommand{\eear}{\end{eqnarray}}
\newcommand{\bean}{\begin{eqnarray*}}
\newcommand{\eean}{\end{eqnarray*}}

\def\AA{{\cal A}}
\def\BB{{\cal B}}

\def\DD{{\cal D}}
\def\EE{{\cal E}}
\def\FF{{\cal F}}

\def\II{{\cal I}}

\def\OO{{\cal O}}

\def\RR{{\cal R}}

\def\eps{\varepsilon}
\def\wlim{\hbox{w-lim} \,}
\def\rlim{r\hbox{-}lim \,}
\def\rlimsup{r\hbox{-}limsup \,}
\def\rliminf{r\hbox{-}liminf \,}
\def\rto{\mathop{\rightharpoonup}^r}
\def\wwto{\displaystyle{\mathop{\rightharpoonup}^{ww}}}
\def\bto{\mathop{\rightharpoonup}^{Ch}}
\def\bto{\mathop{\rightharpoonup}^{b}}
\def\wto{\rightharpoonup}
\def\gamam{\gamma_{\! _-}}
\def\gamap{\gamma_{\! _+}}
\def\sigx{\sigma_{\! x}}

\begin{document}

\title{Kinetic equations with Maxwell boundary conditions} 

\date{}

\author{S. {\sc Mischler}$^1$}

\footnotetext[1]{CEREMADE, Universit\'e Paris IX-Dauphine,
Place du Mar\'echal de Lattre de Tassigny, 75775 Paris, France. 
E-mail: \texttt{mischler@ceremade.dauphine.fr}}

\maketitle

\begin{abstract} 
We prove global stability results of {\sl DiPerna-Lions} renormalized
solutions for the initial boundary value problem associated to some kinetic equations, from which existence results classically follow. 
The (possibly nonlinear) boundary conditions are completely or partially diffuse, which includes the so-called
Maxwell boundary conditions, and we prove that it is realized (it is not only a boundary inequality condition as it has been established in previous works). We are able to deal with Boltzmann, Vlasov-Poisson and Fokker-Planck type models.  The proofs use some trace theorems of the kind previously
introduced by the author for the Vlasov equations, new results concerning weak-weak convergence (the
renormalized convergence and the biting $L^1$-weak convergence), as well as the Darroz\`es-Guiraud information
in a crucial way.
\end{abstract}

\bigskip\bigskip

\centerline{\bf \'Equations cin\'etiques avec conditions aux limites de Maxwell}

\bigskip\noindent{\sl R\'esum\'e - }
Nous montrons la stabilit\'e des solutions renormalis\'ees au sens de {\sl DiPerna-Lions} pour des
\'equations cin\'etiques avec conditions initiale et aux limites.
La condition aux limites (qui peut \^etre non lin\'eaire) est partiellement diffuse et est r\'ealis\'ee
(c'est-\`a-dire qu'elle n'est pas relax\'ee). Les techniques que nous introduisont sont illustr\'ees sur l'\'equation de
Fokker-Planck-Boltzmann et le syst\`eme de Vlasov-Poisson-Fokker-Planck ainsi que pour des
conditions aux limites lin\'eaires sur l'\'equation de Boltzmann et le syst\`eme de Vlasov-Poisson. 
Les d\'emonstrations utilisent des th\'eor\`emes de trace du type de ceux introduits par l'auteur
pour les \'equations de Vlasov, des r\'esultats d'Analyse Fonctionnelle sur les convergences
faible-faible (la convergence renormalis\'ee et la convergence au sens du {\it biting Lemma}),
ainsi que l'information de Darro\`es-Guiraud d'une mani\`ere essentielle.

\bigskip

\noindent
\textbf{Mathematics Subject Classification (2000)}: 76P05 Rarefied gas
flows, Boltzmann equation [See also 82B40, 82C40, 82D05].

\smallskip

\noindent
\textbf{Keywords}: Vlasov-Poisson, Boltzmann and Fokker-Planck equations, Maxwell or diffuse
reflection, nonlinear gas-surface reflection laws, Darroz\`es-Guiraud information, trace Theorems, renormalized
convergence, biting Lemma, Dunford-Pettis Lemma..

\vspace{0.3cm}

\tableofcontents

\section{Introduction and main results}
\setcounter{equation}{0}
\setcounter{theo}{0}

\smallskip
Let $\Omega$ be an open and bounded subset of $\R^N$ and set $\OO = \Omega \times \R^N$.
We consider a gas confined in $\Omega \subset \R^N$. The state of the gas is given by the
distribution function $f = f(t,x,{v}) \ge 0$ of particles, which at time $t \ge 0$ and at position
$x \in \Omega$, move with the velocity ${v} \in \R^N$. 
The evolution of $f$ is governed by a kinetic equation written in the domain $(0,\infty) \times \OO$ and it is complemented with a boundary condition that we describe now. 

We assume that the boundary $\partial\Omega$ is sufficiently smooth. The regularity that we need is that there
exists a vector field $n \in W^{2,\infty}(\Omega; \R^N)$ such that $n(x)$ coincides with
the outward unit normal vector at $x \in \partial\Omega$. We then define $\Sigma_\pm^x := \{ {v} \in \R^N; \pm \, {v} \cdot n(x) > 0 \}$ the sets of outgoing ($\Sigma_+^x$) and incoming ($\Sigma_-^x$) velocities at point $x \in \partial\Omega$ as well as $\Sigma = \partial\Omega \times \R^N$ and 
$$
\Sigma_\pm = \{ (x,{v}) \in \Sigma; \pm n(x) \cdot {v} > 0 \} = \{Ê(x,{v}); \, x \in \partial\Omega, \, {v} \in \Sigma^x_\pm \}. 
$$
We also denote by $d\sigma_{\! x}$ the Lebesgue surface measure on $\partial\Omega$
and by $d\lambda_k$ the measure on $(0,\infty) \times \Sigma$ defined by $d\lambda_k = |n(x) \cdot
{v} |^k \, dt d\sigma_{\! x} dv$, $k=1$ or $2$.

The boundary condition takes into account how the particles are reflected by the wall and thus takes the form of a balance between the values of the trace $\gamma f$ of $f$ on the outgoing and incoming velocities subsets of the boundary:
\beqn\label{BdaryCond}
(\gamma_- f )(t,x,{v}) =  \RR_x( \gamma_+ f (t,x,.) ) ({v}) \,\, \hbox{ on }  \,\,(0,\infty) \times \Sigma_-, 
\eeqn
where $\gamma_\pm f := {\bf 1}_{(0,\infty) \times \Sigma_\pm} \, \gamma f$. The reflection operator is time independent, local in position but can be local or nonlocal in the velocity variable. 
In order to describe the interaction between particles and wall by the mean of the reflection operator $\RR$, J. C. Maxwell \cite{52} proposed in 1879 the following phenomenological law by splitting the reflection operator into a local reflection operator  and a  diffuse  (also denominated as {\it Maxwell}) reflection operator (which is  nonlocal in the velocity variable):
\beqn\label{RLD}
\RR =  (1-\alpha) \, L + \alpha \, D .
\eeqn
Here $\alpha \in (0,1]$ is a constant, called the {\it accommodation coefficient}. The local
reflection operator $L$ is defined by 
$$
(L_x \, \phi) \, ({v}) = \phi (R_x \, {v}),
$$
with $R_x \, {v} = -{v}$ (inverse reflection) or $R_x \, {v} = {v} - 2 \, ({v} \cdot n(x)) \, n(x)$ (specular reflection). 
The diffuse reflection operator $D = (D_x)_{x\in\partial\Omega}$ according to the Maxwellian profile $M$  with temperature (of the wall) $\Theta > 0$  is defined at the boundary point $x \in \partial\Omega$ for any measurable function $\phi$ on $\Sigma^x_+$ by 
$$
(D_x \, \phi) ({v}) = M({v}) \,\, \tilde \phi (x), 
$$
where the normalized Maxwellian $M$ is 
\beqn\label{defMxi}
M({v}) = (2 \, \pi)^{1-N \over 2} \, \Theta^{- {N+1 \over 2} } \, e^{- {|{v}|^2 \over 2 \,
\Theta} }, 
\eeqn
and the out-coming flux of mass of particles $\tilde \phi(x)$  is 
\beqn\label{deftildephi}
\tilde \phi(x) = \int_{{v}' \cdot n(x) > 0} \phi({v}') \, {v}' \cdot n(x) \, d{v}' = \int_{\Sigma_+^x} {\phi \over M}Ê\, d\mu_x.
\eeqn
It is worth emphasizing that the normalization condition (\ref{defMxi}) is made in order that the measure $d\mu_x ({v}) := M({v}) \, |n(x) \cdot {v}| \, d{v}$ is a probability measure on $\Sigma^x_\pm$ for any $x \in \partial\Omega$. Moreover,  for any measurable function $\phi$ on $\Sigma^x_+$ there holds
\beqn\label{OpConsMass}\qquad
\int_{\Sigma_-^x}Ê\!\!\RR_x \phi \, |n(x) \cdot {v}| \, d{v} 
= \int_{\Sigma_-^x}Ê\!\!L_x \phi \, |n(x) \cdot {v}| \, d{v} 
= \int_{\Sigma_-^x}Ê \!\!D_x \phi \, |n(x) \cdot {v}| \, d{v} 
= \int_{\Sigma_+^x}Ê\! \phi \, n(x) \cdot {v} \, d{v}, 
\eeqn
which means that all the particles which reach the boundary are reflected (no particle goes out of the domain nor enters in the domain). 

The reflection law (\ref{RLD}) was the only model for the gas/surface interaction that appeared in the literature before the
late 1960s. In order to describe with more accuracy the interaction between molecules and wall,
other models have been proposed in \cite{26}, \cite{27}, \cite{49} where the reflection operator $\RR$ is a general
integral operator satisfying the so-called non-negative, normalization and reciprocity
conditions, see \cite{30} and Remark~\ref{RGal}.  We do not know whether our analysis can be adapted to such a general kernel. However, the boundary condition can be generalized in an other direction, see \cite{31}, \cite{12}, and we will sometimes assume that the following nonlinear boundary condition holds
\beqn
\label{NLBC}
\RR \, \phi =  (1-\tilde\alpha) \, L \, \phi +  \tilde\alpha \, D \, \phi,
\qquad   \tilde\alpha = \alpha (\tilde \phi),
\eeqn
where $\alpha : \R_+ \to \R_+$ is a continuous function which satisfies $Ê0 < \bar\alpha \le \alpha(s) \le 1$ for
any $s \in \R_+$.

\smallskip
In the domain, the evolution of $f$ is governed by a kinetic equation
\beqn\label{KEq}
{\partial f \over \partial t}  +  {v} \cdot \nabla_x  f = \II (f)
\quad\hbox{ in }\quad (0,\infty) \times \OO,
\eeqn
where $\II(f)$ models the interactions of particles each one with each other and with the environment. Typically, it may be a combination of the quadratic Boltzmann collision operator (describing the collision interactions of particles by binary elastic shocks), the Vlasov-Poisson operator 
(describing the fact that particles interact by the way of the two-body long range Coulomb force) or the Fokker-Planck operator (which takes into account the fact that particles are submitted to a heat bath). More precisely, for the nonlinear boundary condition (\ref{NLBC}) we are able to deal with Fokker-Planck type equations, in particular the Fokker-Planck-Boltzmann equation (FPB in short) and the Vlasov-Poisson-Fokker-Planck system (VPFP in short), while for a constant accommodation coefficient we are able to deal with Vlasov type equations  such as the Boltzmann equation and the Vlasov-Poisson system (VP in short). We refer to section~\ref{sectionModels} where these models are presented. It is worth mentioning that
the method presented in this paper seems to fail for the Vlasov-Maxwell system. 

Finally, we complement these equations with a given initial condition
\beqn\label{CIf0}
f(0,.) = f_{in} \ge 0 \ \hbox{ on } \OO,
\eeqn
which satisfies the natural physical bounds of finite mass, energy and entropy
\beqn\label{bornef0}
\int \!\! \int_\OO f_{in} \, ( 1 + |{v}|^2 + |\log \, f_{in}|) \, dx d{v} =: C_0 < \infty.
\eeqn

We begin with a general existence result that we state deliberately in an imprecise way and we refer to section~\ref{sectionModels} (and Theorem~\ref{Th:ModStab}) for a more precise statement. 

\begin{theo}\label{ThExistence} Consider the initial boundary value problem (\ref{BdaryCond})-(\ref{KEq})-(\ref{CIf0}) associated to the FPB equation or the VPFP system with possibly mass flux depending accommodation coefficient (\ref{NLBC}) or the boundary value problem associated to the Boltzmann equation or the VP system with constant accommodation coefficient (\ref{RLD}).
For any non-negative initial datum $f_{in}$ with finite mass, energy and entropy ((\ref{bornef0}) holds) there exists at least one (renormalized) solution $f \in C([0,\infty);L^1(\OO))$ with finite mass, energy and entropy to  the kinetic equation (\ref{KEq}) associated to the initial datum $f_{in}$ and such that the trace function $\gamma f$ fulfills  the boundary condition (\ref{BdaryCond}).
%(it is not a boundary inequality condition). 
\end{theo}

The Boltzmann equation and the FPB equation for initial data satisfying the natural bound (\ref{bornef0}) was
first studied by R. DiPerna and P.-L. Lions \cite{35,37,38} who proved stability and existence results
for  global renormalized solutions in the case of the whole space ($\Omega = \R^N$).
Afterwards, the corresponding boundary value problem with reflection boundary conditions (\ref{BdaryCond}) and constant accommodation coefficient  has been extensively studied in the case of the Boltzmann model \cite{45}, \cite{4}, \cite{5}, \cite{6}, \cite{7},
\cite{28}, \cite{41}, \cite{46}, \cite{29}, \cite{54}. It has been proved, in the partial absorption case
$\gamam f = \theta \, \RR \gamap f $ with $\theta \in[0,1)$ and in the completely local reflection
case (i.e. (\ref{BdaryCond}) holds with $\alpha \equiv 0$), that there exists a global renormalized solution.
But in the most interesting physical case (when $\theta \equiv 1$ and $\alpha \in (0,1]$), it has
only been proved in \cite{6} that the following boundary inequality condition
\beqn\label{bdConRelax}
\gamma_- f \ge \RR (\gamma_+ f) \quad \hbox{ on } (0,\infty) \times \Sigma_-
\eeqn 
holds, instead of the boundary equality condition (\ref{BdaryCond}). However, it is worth mentioning that if  the renormalized solution built in \cite{6} is in fact a solution to the Boltzmann equation in the sense of distributions, then that solution satisfies the boundary equality condition (\ref{BdaryCond}) (a result that one deduces thanks to the Green formula by gathering the fact that the solution is mass preserving and the fact that the solution already satisfies the boundary inequality condition (\ref{bdConRelax})). Also, the Boltzmann equation with nonlinear boundary conditions has been treated in the setting of a strong but non global solution framework in \cite{42}.

\smallskip
With regard to existence results for the initial value problem for the VPFP system set in the
whole space, we refer to \cite{14}, \cite{15}, \cite{16}, \cite{19}, \cite{22}, \cite{23}, \cite{34}, \cite{57}, \cite{64}, \cite{24}, \cite{59}
as well as \cite{20} for physical motivations. The initial boundary value problem has been addressed in \cite{13}, \cite{21}. We also
refer to \cite{3}, \cite{11}, \cite{44}, \cite{56}, \cite{66} for the initial boundary value problem for the VP system
and to \cite{56} for the corresponding stationary problem. We emphasize that in all these
works only local reflection or prescribed  incoming data are treated, and to our knowledge, there
is no result concerning the diffuse boundary condition for the VP system or for the VPFP system.

\smallskip
We also mention that there is a great deal of information for the boundary value problem in an
abstract setting in \cite{65}, \cite{43} with possibly nonlinear boundary conditions \cite{10}, \cite{55}. 

\smallskip
In short, the present work improves the already known existence results for kinetic equations with diffusive boundary reflection into three directions. 

\smallskip
$\bullet$ On the one hand, we prove that (\ref{BdaryCond}) is fulfilled, while only the boundary inequality condition (\ref{bdConRelax}) was previously  established.

\smallskip
$\bullet$ On the other hand, we are able to consider a large class of kinetic models (including Vlasov-Poisson term) while only the Boltzmann equation (or linear equations) could be handled with earlier techniques. 
%{\bf (earlier these earlier techniques)???} 

\smallskip
$\bullet$ Finally, we are able to handle some nonlinear boundary condition in the case of Fokker-Planck type equation.

\smallskip
We do not present the proof of Theorem~\ref{ThExistence} (nor the proof of its accurate version 
Theorem~\ref{Th:ModStab}) because it classically follows from a sequential stability 
or sequential compactness result that we present below and a standard (but tedious) approximation 
procedure, see for instance \cite{54} or the above quoted references. We deliberately  state again the   sequential stability  result in  an imprecise way, referring to section~\ref{sectionModels} for a more accurate version. 

\begin{theo}\label{ThStab} Consider the initial boundary value problem (\ref{BdaryCond})-(\ref{KEq})-(\ref{CIf0}) associated to the FPB equation or the VPFP system with possibly mass flux depending constant accommodation coefficient (\ref{NLBC}) or the boundary value problem  associated to the Boltzmann equation or the VP system with constant accommodation coefficient (\ref{RLD}).
Let then $(f_n)$ be a sequence of (renormalized) solutions to that equation and assume that $(f_n)$ and the trace sequence $(\gamma f_n)$ satisfy the natural physical a priori bounds (to be specified for each model). If $f_n(0,.)$ converges to $f_{in}$ weakly in $L^1(\OO)$ then, up to the extraction of a subsequence, $f^n$ converges (at least) weakly in $L^1([0,T]Ê\times \OO)$ for all $T \in (0,\infty)$ to a (renormalized) solution $f$ to the kinetic equation (\ref{KEq}) with initial value $f_{in}$. 
Furthermore, for any $\eps > 0$ and $T > 0$ there exists a measurable set $A \subset (0,T) \times \partial\Omega$ such that 
meas$\, ((0,T) \times \partial\Omega \, \backslash \, A) < \eps$ and 
\beqn\label{cvgeChaWeakxi}
\gamap f_n \wto\gamap f \quad\hbox{weakly  in}\quad L^1(A \times \R^N, d\lambda_1),
\eeqn
the convergence being strong in the case of the Fokker-Planck type equations. 
As a consequence we can pass to the limit in the reflection boundary condition  (\ref{BdaryCond})-(\ref{RLD}) (and (\ref{BdaryCond})-(\ref{NLBC}) in the case of Fokker-Planck type equations), so that the reflection boundary condition (\ref{BdaryCond}) is fulfilled. 
\end{theo}

\smallskip
Let us briefly explain the main steps and difficulties in the proof of the stability result.

\smallskip
$\bullet$ The first step consists in collecting the physical estimates available on the solution $f$ to the equation (\ref{BdaryCond})-(\ref{KEq})-(\ref{CIf0}) and on its trace $\gamma f$. In the interior of the domain the a priori bounds satisfied by $f$  strongly depend on the model considered but  they are the same than those available in the case of the whole space. In general, for the trace, we are only able to prove that 
\beqn\label{DGinfo}
\forall \, T \qquad \int_0^T \!\! \int_{\partial\Omega} \EE\bigr( {\gamma_+ f \over M} \bigr) \,  d\sigx dt \le C_T,
\eeqn
with $C_T$ only depending on $C_0$ and $T$, where the functional  $\EE = \EE_x$ is the {\it Darroz\`es-Guiraud information} defined by 
\beqn\label{defDGI}
\EE\bigl( \phi \bigr) :=  \int_{ \Sigma_+^x } h ( \phi) \, d\mu_x - h \left(  \int_{ \Sigma_+^x } \phi \, d\mu_x \right), \qquad
h(s) = s \, \log s,
\eeqn
and where we recall that $d\mu_x({v}) :=  M({v}) \, |n(x) \cdot {v}| \, d{v}$ is a probability measure on $\Sigma_+^x$ so that $\EE(\phi) \ge 0$ thanks to the Jensen inequality. 
Let us emphasize that additionally to the a priori bound of the {\sl Daroz\`es-Guiraud
information} (\ref{DGinfo}), we can prove an $L^1$ a priori bound in the case of the Boltzmann
equation (and of the FPB equation) and only an $L^{1/2}$ a priori  (but also a posteriori) bound in the case of the VP system (and the FPVP system): in both cases, we do not have any a priori information on the trace which guaranties  uniform local equiintegrability on the trace functions of a sequence of solutions.  The main difficulty is thus the lack of a good a priori bound on the trace. 

\smallskip
$\bullet$ The next step consists in specifying the sense of the equations. The physical a priori estimates on $f$ make possible to give a sense to (\ref{KEq}) in a renormalized sense as introduced by DiPerna and Lions. What is then the meaning of the trace $\gamma f$ of $f$?  That  so-called trace problem has been studied in \cite{9}, \cite{32}, \cite{2}, \cite{62}, \cite{43}, \cite{18} for the Vlasov equation with a Lipschitz force field and extended to the Vlasov-Fokker-Planck equation in \cite{21}. 
In the case of the VP and the VPFP systems, the a priori estimate on the force field does not guarantee Lipschitz regularity but only Sobolev regularity. A trace theory has been developed in \cite{53,54} for the (possibly renormalized) solutions of the Vlasov equation with a force field in Sobolev space that we extend here to the solutions of the
Vlasov-Fokker-Planck equation. The trace of a solution is here defined by a Green formula written on the
renormalized equation.

\smallskip
$\bullet$ In a last step, we have to pass to the limit in a sequence of solutions which satisfy the ``natural physical bounds".  For the equation satisfied by $f$ in the interior of the domain, the proofs have been done already by  DiPerna-Lions \cite{34,35,37} and Lions \cite{50}, and nothing has to be changed.  The main difficulty solved here is to handle the boundary condition which is made up of two equations:

\smallskip
(1) the renormalized Green formula which links together the solution $f$ in the interior of the domain with its trace function $\gamma f$; 

\smallskip

(2) the boundary equality condition (\ref{BdaryCond}) which connects together the incoming velocity particles density $\gamma_- f$ with the outgoing velocity particles density $\gamma_+ f$.

\smallskip\noindent
Let us emphasize that using only the $L^1$ boundedness information (as it is available for the Boltzmann equation for instance) on a sequence $(\gamma f_n)$ of the trace of solutions to a kinetic equation satisfying the boundary condition (\ref{BdaryCond}) it is only possible to prove the boundary inequality condition (\ref{bdConRelax}). Indeed, on the one hand as in \cite{6} we may use that, up to the extraction of a subsequence, $f_n \wto f$ weakly in $L^1$ and $\gamma_\pm f_n \wto \eta_\pm$ in the weak sense of measures for some measures $\eta_\pm \ge 0$. Then the limit boundary densities $\eta_\pm$ fulfill the boundary equality condition (\ref{BdaryCond}), $\eta_- = \RR\, \eta_+$, whereas they are not the trace functions associated to $f$ but they are their regular parts with respect to the lebesgue measure: $\gamma f_\pm = {d \eta_\pm \over d\lambda_1}$. Putting together these two informations yields to the boundary inequality condition (\ref{bdConRelax}). On the other hand, as in \cite{41}, we may use that, up to the extraction of a subsequence, $f_n \wto f$ weakly in $L^1$ and $\gamma_\pm f_n \wto g_\pm$ in the biting $L^1$-weak sense (see below) for some measurable functions $g_\pm \ge 0$. Then the limit boundary densities $g_\pm$ are the trace functions associated to $f$, $g_\pm = \gamma_\pm f$, whereas the reflection operator is only l.s.c. with respect to the biting $L^1$-weak convergence, $\RR \, g_+ \le \liminf \RR \, \gamma_+ f_n$. Again, these two informations only imply the boundary inequality condition (\ref{bdConRelax}). 

\smallskip
$\bullet$
In this paper, we prove some $L^1$-weak ($L^1$-strong in the case of FP models) convergence in the velocity variable for the sequence $(\gamma_+ f_n)$ (as stated in Theorem~\ref{ThStab}) which is strong enough to conclude. Our proof is based on the use of notions of weak-weak convergences, namely the renormalized convergence (r-convergence) and the biting $L^1$-weak convergence (b-convergence).  We say weak-weak convergences in order to express the fact that they are extremely weak sense of convergence (weaker, for instance, to the $L^1$-weak convergence and to the a.e. convergence) and which are not furthermore associated to any topological structure, see Proposition~\ref{prop:Appendix}. 
On the one hand, thanks to the trace theory, we prove that the sequence of trace functions $(\gamma f_n)$ r-converges to $\gamma f$ (as well as a.e. in the case of  FP models). Next, thanks to some additional $L^1$ a priori bounds, or because the r-convergence is  almost equivalent to the b-convergence when the limit function belongs to $L^0$, we deduce that $\widetilde{\gamma_+ f_n}$ b-converges to $\widetilde{\gamma_+ f}$. Finally, that information and the boundedness of the Darroz\`es-Guiraud information leads to (\ref{cvgeChaWeakxi}). 

\smallskip
Let us now briefly outline the contents of the paper. In section~\ref{sec:FT}, we consider the free transport equation for which we apply the above strategy. We present for this very simple case the different tools (renormalized and biting $L^1$-weak convergence, trace theory and Darroz\`es-Guiraud information), we state a first  velocity $L^1$-weak compactness result and then we prove the corresponding version of the  stability Theorem~\ref{ThStab}. 
In Section~\ref{sec:WW}, we develop the notion of renormalized convergence in a more general framework and we prove some more accurate version of biting $L^1$-weak convergence and velocity $L^1$-weak compactness.
In Section~\ref{sec:Trace}, we present the trace theory for the Vlasov-Fokker-Planck equation with Sobolev regularity on the force field. 
In Section~\ref{sec:GalStab}, putting together the results from Section~\ref{sec:WW} and Section~\ref{sec:Trace}, we establish the renormalized convergence and the almost everywhere convergence of trace functions sequences. In Section~\ref{sec:Models} we present the models and we establish the main stability (up to the boundary) results. Finally, in the Appendix, we come back to the notion of renormalized convergence for which we give several relevant examples and counterexamples.

%%%%%%%%%%%%%%%%%%%%%%%%%%%%%%%%%%%%%%%%%%%%%%%%
\section{An illuminating example: the free transport equation.} \label{sec:FT}
\setcounter{equation}{0}
\setcounter{theo}{0}
%%%%%%%%%%%%%%%%%%%%%%%%%%%%%%%%%%%%%%%%%%%%%%%%

\smallskip
In this section we assume that $f$ is governed by the free transport equation
\beqn\label{FTeq}
{\partial f \over \partial t}  +  {v} \cdot \nabla_x  f =  0 \quad \hbox{ in }\
(0,\infty) \times \OO,
\eeqn
complemented with the initial condition (\ref{CIf0}) and the boundary reflection condition (\ref{BdaryCond})
with constant restitution coefficient $\alpha \in (0,1]$. 
Our aim is  to adapt the {\it DiPerna-Lions stability theory} to that simple boundary value problem. We follow the strategy 
expounded in the introduction. We first collect the a priori bounds satisfied by a solution to the boundary value problem (\ref{FTeq})-(\ref{BdaryCond})-(\ref{CIf0}) with initial datum satisfying (\ref{bornef0}). We next present some general functional analysis tools which  roughly speaking make  possible to deduce the $L^1$ weak convergence in the ${v}$ variable of a sequence which is uniformly bounded in $L^1$ and for which the associated {\it Darroz\`es-Guiraux information} is uniformly bounded. We finally state and prove the stability result associated to  the boundary value problem (\ref{FTeq})-(\ref{BdaryCond})-(\ref{CIf0}). 

\begin{rem}\label{RemBoltz}
It is worth mentioning that  the proof of the corresponding stability result for the Boltzmann equation is essentially the same as for the free transport equation. We refer to section~6 where that model is handled. However, the reader who is only interested in the Boltzmann model may easily adapt the proof below with the arguments introduced in \cite{53} (it will be more elementary than the proof presented in Section~3 to Section~6 which is made in order to also deal with a Vlasov-Poisson term and/or with a Fokker-Planck term). 
\end{rem}

%\begin{rem}\label{referee}
%Dans le th 2.9.  
%On veut etablir un theoreme de "compacite/stabilite": etant donnee une suite de solutions qui satisfont les bornes naturelles uniformement, alors on peut en extraire une sous-suite qui converge vers une limite qui est encore une solution. 
%1) Le premier point est donc de determiner quelles sont les bornes;
%Ensuite le problme que l'on se pose est de passer a la limite "au bord", le passage a l'interieur ayant deja ete traite ailleurs. 
%2) On extrait de 1) une convergence
%3) il faut passer ensuite a la limite dans LES DEUX EQUATIONS satisfaites par $\gamma f$: 
%  3a - la formule de Green qui relie $\gamma f$ a $f$;
%  3b - l'equation de reflexion.
%Le pb est qu'en utilisant uniquement la borne $L^1$, comme cela a ete fait anterieurement, on ne sait en dŽduire qu'une convergence au sens biting, 
%et donc on obtient au mieux que $\eta_- \ge R \eta_+$. La remarque fondamentale est donc d'observer que la borne $L^1$ $+$ la borne sur l'information de DG donne une convergence au sens biting en $t,x$ et $L^1$ faible en $v$. On obtient ainsi Pour passer a la limite dans la formule de Green, 

%\end{rem}

%%%%%%%%%%%%%%%%%%%%%%%%%%%%%%%%%%%%%%%%%%%%%%%

\subsection{A priori bounds.}

\begin{lem}\label{lem:FreeTApBound}
For any non-negative initial datum $f_{in}$ such that (\ref{bornef0}) holds
and any time $T \in (0,\infty)$ there exists a constant $C_T$ (only depending on $C_0$ and $T$) such that any sufficiently regular and decreasing at the infinity solution $f$ to the initial boundary value problem (\ref{FTeq})-(\ref{BdaryCond})-(\ref{CIf0}) satisfies %(at least formally) 
\beqn\label{FTbd1}
\sup_{[0,T]}Ê\int \!\! \int_\OO f \, ( 1 + |{v}|^2 + |\log \, f|) \, dx d{v}
+ \alpha \int_0^T \!\int_{\partial\Omega} \EE \left( {\gamma_+ f \over M}Ê\right) \, d\sigx dt  \le C_T,
\eeqn
where $\EE$ is defined in (\ref{defDGI}), and 
\beqn\label{FTbd2}
 \alpha \int_0^T \int\!\!\int_\Sigma \gamma f \, (1 + |{v}|^2) \, |n(x) \cdot {v}| \, d{v} d\sigma_{\! x} dt
 \le \, C_T.
\eeqn
\end{lem}

\medskip\noindent
{\sl Proof of Lemma~\ref{lem:FreeTApBound}. } We consider a solution $f$ of  (\ref{BdaryCond})-(\ref{FTeq})-(\ref{CIf0}), which is sufficiently regular and decreasing at the infinity in such a way that all the integrations by parts in our arguments are legitimate.

\noindent
{\sl Step 1. Mass conservation. }ÊIntegrating the free transport equation (\ref{FTeq}) over $x,{v}$, using the Green formula and the  identity (\ref{OpConsMass}), we obtain the mass conservation
$$
\forall \, t \ge 0 \qquad \int\!\!\int_\OO f(t,.) \, d{v} dx = \int\!\!\int_\OO f_{in} \, d{v} dx.
$$ 

\medskip\noindent
{\sl Step 2. Relative entropy. } Multiplying  the free transport  equation (\ref{FTeq}) by $h'(f/M)$, with $h(s) = s \log s$, and integrating it over $x,{v}$, we have 
\beqn\label{FTbd3}
{d \over dt}  \int \!\! \int_\OO h(f/M) \, M \, d{v} dx =  \int \!\! \int_\Sigma 
h(\gamma f/M) \, M \, {v} \cdot n(x) \,\, d{v} d\sigx.
\eeqn
The Darroz\`es-Guiraud inequality states that the entropy boundary flux at the right hand side of equation (\ref{FTbd3}) is non-negative. That is a straightforward consequence of the Jensen inequality taking advantage that $d\mu_x(v) =  M \, |{v} \cdot n(x)| \,\, d{v}$ is a probability measure. We present now the proof of an accurate version of the Darroz\`es-Guiraud inequality which make precise how much that term is non-negative. 
From  the boundary reflection condition (\ref{BdaryCond}), the convexity of $h$ and the expression (\ref{RLD}) of the reflection operator, we have 
\bear\label{FTbd4}
&&\int_{\R^N} \!\! h(\gamma f/M) \, d\mu_x({v})  
= \int_{ \Sigma^+_x} h ( \gamma_+ \, f/M) \, d\mu_x({v}) - \int_{ \Sigma^-_x} h ( \RR \gamma_+ f/M) \, d\mu_x({v}) \\
\nonumber
&&\qquad \ge \alpha \left\{Ê\int_{ \Sigma^+_x} h ( \gamma_+ \, f/M) \, d\mu_x({v}) - \int_{ \Sigma^-_x} h ( D \gamma_+ f/M) \, d\mu_x({v}) \right\} \\ \nonumber
&&\qquad+ (1-\alpha) \left\{Ê\int_{ \Sigma^+_x} h ( \gamma_+ \, f/M) \, d\mu_x({v}) - \int_{ \Sigma^-_x} h ( L \gamma_+ f/M) \, d\mu_x({v}) \right\} \\ \nonumber
&&\qquad = \alpha \left\{ \int_{ \Sigma^+_x} Ê h ( \gamma_+ \, f/M) Ê\, d\mu_x({v})  - h (\widetilde{\gamap f }) \right\} =
\alpha \, \EE_{x}({\gamma_+ f \over M}),
\eear
where we have performed the change of variables $L_x : {v} \mapsto R_x {v}$ in the second term with $jac \, L_x = 1$, so that this term vanishes, and where the  Darroz\`es-Guirraud information functional $\EE_x$ is defined in (\ref{defDGI}) and $\widetilde{\gamap f }$ is defined in (\ref{deftildephi}). Gathering (\ref{FTbd3}) and (\ref{FTbd4}), we get 
$$
{d \over dt}  \int \!\!\! \int_\OO h(f/M) \, M \, d{v} dx + \alpha \, \int_{\partial\Omega} \EE_x ( \gamma_+ \, f) \, d\sigx  
 \le 0.
$$
Finally, using the elementary estimates, that one can find in \cite{51} for instance, 
\beqn\label{BddH1}
\int_{\R^N} f \, \Bigl( {|{v}|^2 \over 4 \, \Theta} + | \log f| \Bigr) \, d{v} 
\le C_M + \int_{\R^N}h(f/M) \, M \, d{v},
\eeqn
and 
\beqn\label{BddH2}
\int_{\R^N}h(f_{in}/M) \, M \, d{v} 
\le \int_{\R^N} f_{in}\, \Bigl( {|{v}|^2 \over 4 \, \Theta} + | \log f_{in}| \Bigr) \, d{v} + C_M,
\eeqn
for some constant $C_M \in (0,\infty)$, we obtain that (\ref{FTbd1}) holds. 

\medskip\noindent
{\sl Step 3. Additional $L^1$ estimates. } 
For the sake of completeness we sketch the proof of the $L^1$ a priori bound (\ref{FTbd2}) already established in \cite{6,54}.
We multiply the free transport equation (\ref{FTeq}) by $n(x) \cdot {v}$ and we integrate it over all variables, to get
$$
\int_0^T \!\!\! \int \!\!\! \int_\Sigma \gamma f \, (n(x) \cdot {v})^2 \, d{v} d\sigx dt
= \left[ \int \!\!\! \int_\OO f \, n(x) \cdot {v} \, d{v} dx \right]_T^0 + 
\int_0^T \!\!\! \int \!\!\! \int_\OO f \, {v} \cdot \nabla_x n(x)  {v} \, d{v} dx dt,
$$ 
so that, thanks to (\ref{FTbd1}) and because $n \in W^{1,\infty}(\Omega)$, 
\beqn\label{FTbd5}
\int_0^T \!\!\! \int \!\!\! \int_\Sigma \gamma f \, (n(x) \cdot {v})^2 \, d{v} d\sigx dt \le C_T.
\eeqn 
We then remark that  for the constant $C_1 :=  \|ÊM(v) \, (n(x) \cdot v)^2 \|_{L^1(\Sigma^x_-)}^{-1}$ we have 
\bear\label{FTbd6}
\widetilde{\gamap f} 
= C_1  \int_{\Sigma_-^x} 
M({v}) \, \widetilde{\gamap f}  \, (n(x) \cdot {v})^2 \, d{v} 
= C_1 \int_{\Sigma_-^x} 
\gamam f \, (n(x) \cdot {v})^2 \, d{v}, 
\eear
and that  for the constant $C_2 :=  \|ÊM(v) \, (1+|v|^2) \, |n(x) \cdot v| \|_{L^1(\Sigma^x_-)}$ we have 
\beqn\label{FTbd7}
\int_{\Sigma_-^x} \gamam f \, (1+ |{v}|)^2 \, |n(x) \cdot {v}| \, d{v}
= \int_{\Sigma_-^x} M({v}) \, \widetilde{\gamap f} \,(1+  |{v}|^2) \, |n(x) \cdot {v}| \, d{v}
= C_2 \, \widetilde{\gamap f}. 
\eeqn
Finally, we come back to the equation (\ref{FTeq}) that we multiply by $|{v}|^2$ and that we integrate in all variables. We obtain
\bear\label{FTbd8}
\int \!\!\! \int_\OO f(T,.) \, |{v}|^2 \, d{v} dx 
&+&\alpha\int_0^T \!\!\! \int \!\!\! \int_{\Sigma_+} \gamap f \, |{v}|^2 \, n(x) \cdot {v} \, d{v} d\sigx dt \\ \nonumber
&= &\int \!\!\! \int_\OO f_{in} \, |{v}|^2 \, d{v} dx
+\alpha\int_0^T \!\!\! \int \!\!\! \int_{\Sigma_-} \gamam f \, |{v}|^2 \, |n(x) \cdot {v}| \, d{v} d\sigx dt. 
\eear
Estimate (\ref{FTbd2}) follows gathering (\ref{FTbd6}), (\ref{FTbd5}) and (\ref{FTbd7}), (\ref{FTbd8}). \qed

%%%%%%%%%%%%%%%%%%%%%%%%%%%%%%%%%%%%%%%%%%%%%%%

\subsection{Biting  $L^1$-weak convergence and $L^1$-weak compactness in the velocity variable.}

In this section we present some functional analysis results which make possible to obtain the $L^1$-weak 
convergence in the ${v}$ variable of a sequence which satisfies a  $L^1$ bound and a uniform bound of its {\sl
Darroz\`es-Guiraud information}. We state the result in some more general setting because we believe that it 
may have its own interest (outside the applications to the trace theory for kinetic equations). 
For that purpose, we introduce a first notion of weak-weak convergence, 
namely the biting $L^1$-weak convergence. It seems to have been introduced by Kadec and Pelzy\'nski \cite{48} and rediscovered and developed in a
$L^1$ and bounded measure framework by Chacon and Rosenthal in the end of the 1970's, see \cite{40},
\cite{17}. Let us first recall the definition of the biting $L^1$-weak convergence that we extend to a ``$L$ framework". 

In the sequel $Y = (Y,\nu)$ stands for a separable and $\sigma$-compact topological space, i.e. $Y = \cup_k
Y_k$ where $(Y_k)$ is an increasing sequence of compact sets, endowed  with its $\sigma$-ring
of Borel sets and with a locally finite Borel measure $\nu$. We denote by $L(Y)$ the space of all measurable
functions $\phi: Y \to \bar\R$ and by $L^0(Y)$ the subset of all measurable and
$\nu$-almost everywhere finite functions. In order to simplify the presentation,
we will be only concerned with non-negative functions of $L$ and $L^0$. Thus, in this section, we also
denote by  $L$ and $L^0$ the cone of non-negative functions in these spaces, and we do not 
specify it anymore.

\begin{defin}\label{defbiting}  We say that a sequence $(\psi_n)$ of $L(Y)$ converges in the biting
$L^1$-weak sense (or b-converges) to $\psi \in L(Y)$, denoted $\displaystyle{ \psi_{n} \bto \psi }$, 
if for every  $k \in \N$ we can find $A_k \subset Y_k$ 
in such a way that $(A_k)$ is increasing,
$\nu(Y_k \backslash A_k) < 1/k$, $\psi_n \in L^1(A_k)$ for all $n$ large enough and
$\psi_n \wto \psi$ weakly in $L^1(A_k)$. In particular, that implies $\psi \in L^0(Y)$.
\end{defin}

The fundamental result concerning the biting $L^1$-weak convergence is the so-called biting Lemma
that we recall now. We refer to \cite{25}, \cite{8}, \cite{17}, \cite{40} and \cite{48} for a proof of this
Lemma.  We also refer to \cite{1} and \cite{33} for other developments related to the biting $L^1$-weak
convergence. Extension of this theory to multi-valued functions has been done by Balder, Castaing,
Valadier and others; we refer to \cite{58} for precise references.

\begin{theo}\label{theobiting} 
{\bf (biting Lemma).}  Let $(\psi_n)$ be a bounded sequence of $L^1(Y)$. 
There exists $\psi \in L^1(Y)$ and a subsequence $(\psi_{n'})$ such that 
$(\psi_{n'})$ b-converges to $\psi$ and $\|Ê\psi \|_{L^1} \le \liminf \|Ê\psi_{n'} \|_{L^1}$. 
\end{theo}

\medskip
Our first result is a kind of intermediate result between the biting Lemma and the Dunford-Pettis
Lemma. More precisely, we prove the $L^1$-weak compactness in the ${v}$ variable for sequences
$(\phi_n)$ which are bounded in $L^1$ and such that the associated {\sl Darroz\`es-Guiraud information} 
is uniformly (in $n$) bounded. It is based on the biting Lemma, the Dunford-Pettis Lemma and a convexity argument.

\begin{theo}\label{theo1xi} 
Consider $j : \R_+ \to \R$ a convex function of class $C^2(0,\infty)$ such that 
$j(s)/s \to +\infty$ when $s \nearrow +\infty$ and such that the application $J$ from $(\R_+)^2$ to
$\R$ defined by $J(s,t) = (j(t) - j(s)) \, (t-s)$ is convex, $\omega$ a non-negative function
of $\R^N$ such that $\omega({v}) \to \infty$ when $|{v}| \to \infty$ and, for any $y \in Y$, a
probability measure $\mu_y$ on $\R^N$. Assume that $(\phi_n)$ is a sequence of non-negative
measurable functions on $Y \times \R^N$ such that
\beqn\label{phiBdd1}
\int_Y \!\! \int_{\R^N} \bigl[ \phi_n(y,{v}) \, (1 +Ê\omega({v})) + \EE(\phi_n(y,.)) \bigr] \,
d\mu_y({v}) \, d\nu(y) \le C_1 < \infty,
\eeqn
where $\EE = \EE_{j,y}$ is the non-negative {\sl Jensen information functional} defined by
$$
\EE(\phi) = \int_{\R^N} j(\phi) \, d\mu_y - j \Bigl( \int_{\R^N} \phi \, d\mu_y \Bigr)
\qquad \hbox{ if } \quad 0 \le \phi \in L^1(\R^N,d\mu_y).
$$
Then, there exists $\phi \in L^1(Y \times \R^N)$ and a subsequence $(\phi_{n'})$ such that
for every $k \in \N$ we can find $A_k \subset Y_k$ in such a way
that $(A_k)$ is increasing, $\nu(Y_k \backslash A_k) < 1/k$ and
$$
\phi_{n'} \rightharpoonup \phi \quad\hbox{weakly in}\quad L^1(A_k \times \R^N;d\nu \, d\mu).
$$
Furthermore, $\EE$ is a convex and weakly $L^1$ l.s.c. functional, and thus
\beqn\label{phiBdd2}
\int_Y \!\! \int_{\R^N} \bigl[ \phi(y,{v}) \, (1 +Ê\omega({v})) + \EE(\phi(y,.)) \bigr] \,
d\mu_y({v}) \, d\nu(y) \le C_1.
\eeqn
\end{theo}

\noindent{\sl Proof of Theorem~\ref{theo1xi}.}
From the bound (\ref{phiBdd1}) and the biting Lemma we know that there exists a subsequence $n'$ such that for
every $k \in \N$ we can find a Borel set $A = A_{k} \subset Y_k$ with $\nu(Y_k \backslash A) <
1/k$ such that
\beqn\label{phiBdd3}
\int_{\R^N} \phi_{n'} \, d\mu_y({v}) \quad\hbox{ weakly converges in } \ L^1(A).
\eeqn

Thanks to (\ref{phiBdd3}), the Dunford Pettis Lemma and the De La Vall\'ee-Poussin uniform integrability criterion  there is a convex function $\Phi = \Phi_k$
such that $\Phi (s) / s \to \infty$ when $s \to \infty$ and 
$$
\int_A \Phi \Bigl( \int_{\R^N} \phi_{n'} \, d\mu_y({v}) \Bigr)  \, d\nu(y) \le
C_2 = C_2(k) < \infty.
$$
Furthermore, we can assume that $\Phi(0) = 0$, $\Phi' = a_m$ in $[m,m+1]$ with 
$j'(s_0) \le a_m \nearrow + \infty$, where $s_0 \in \N^\star$ is such that $j(s_0) \ge 0$ and
$j'(s_0) \ge 0$. 

Then we define $\Psi = \Psi_k$ by $\Psi(s) = j(s)$ for $s \in [0,s_0]$ and by induction on 
$m \in \N$, we consider $t_m$ such that $j'(t_m) = a_m - \Psi'(s_m) + j'(s_m)$ and we set
$s_{m+1} = [t_m] + 1$, $\Psi'' := j''$ on $[s_m,t_m]$ and
$\Psi'' := 0$ on $[t_m,s_{m+1}]$ so that $t_m \ge s_m \ge m$ and $\Psi'(s_{m+1}) \ge
a_m \ge \Psi'(s_m)$. Therefore, we have built a convex function $\Psi$ such that the function
$s \mapsto j(s) - \Psi(s)$ is convex, $\Psi (s) / s \nearrow \infty$ since $\Psi'(s) \nearrow
\infty$, and
$\Psi \le \Phi$ since $\Psi' \le \Phi'$, so that 
\beqn\label{phiBdd4}
\int_A \Psi \Bigl( \int_{\R^N} \phi_{n'} \, d\mu \Bigr) 
\, d\nu \le C_2.
\eeqn
The Jensen inequality, written for the function $s \mapsto j(s) - \Psi(s)$, gives 
$$
\int_{\R^N} \Psi(\phi_{n'}) \, d\mu - 
\Psi \bigl( \int_{\R^N} \phi_{n'} \, d\mu \bigr) 
\le \EE(\phi_{n'}),
\leqno
$$
and combining it with (\ref{phiBdd1}) and (\ref{phiBdd4}) we get
$$
\int \!\!\! \int_{A \times \R^N}  \Psi (\phi_{n'}) \, d\mu \, d\nu \le C_1 + C_2,
\leqno
$$
and thus
\bear\label{Psi+}
\int \!\!\! \int_{A \times \R^N}  \Psi^+ (\phi_{n'}) \, d\mu_y \, d\nu
&&\le C_1 + C_2 + \int \!\!\! \int_{A \times \R^N}  \Psi^- (\phi_{n'}) \, d\mu_y \, d\nu \\ \nonumber
&&\le C_3(k) := C_1 + C_2 + \nu(A) \, \sup j^- < \infty.
\eear
Thanks to estimates (\ref{phiBdd1}), (\ref{Psi+}) and  the Dunford-Pettis Lemma we get that
$(\phi_{n'})$ falls in a relatively weakly compact set of $L^1(A_k \times \R^N)$ for any $k \in \N$.
We conclude, by a diagonal process, that there is a function $\phi \in L^1(Y \times \R^N)$ and a
subsequence $(\phi_{n''})$ which converges to $\phi$ in the sense stated in Theorem~\ref{theo1xi}.

In order to prove that $\EE$ is a convex functional, we begin by assuming that $j \in C^1(\R_+,\R)$,
so that $\EE$ is G\^ateaux differentiable. By definition of the G-differential
\bean
\nabla \EE(\phi) \cdot \psi 
&:=& \lim_{t \to 0} {\EE (\phi + t \, \psi) - \EE(\phi) \over t}Ê\\
&=& \int_{ \R^N } j'(\phi) \, \psi \, d\mu - j'\bigl( \int_{ \R^N } \phi \, d\mu \bigr) 
\int_{ \R^N } \, \psi \, d\mu, 
\eean
for any $0 \le \phi,\psi \in L^\infty(\R^N)$. Therefore, by the Jensen
inequality, we have
$$
\langle \, \nabla \EE(\psi) - \nabla \EE(\phi), \psi - \phi \, \rangle \ = 
\int_{ \R^N } J(\phi,\psi) \, d\mu - J \Bigl(\int_{ \R^N } \phi \, d\mu,\int_{ \R^N } \psi \,
d\mu\Bigr) \ge 0,
$$
so that $\nabla \EE$ is monotone and thus $\EE$ is convex on $L^\infty(\R^N)$: for any $0 \le
\phi,\psi
\in L^\infty(\R^N)$ and any $t \in (0,1)$ 
\beqn\label{2.17}
\EE(\phi +(1-t) \, \psi) \le t \, \EE(\phi) + (1-t) \, \EE(\psi).
\eeqn
When $j \notin C^1(\R_+,\R)$ we define, for any $\eps > 0$, the function $j_\eps(s) = j(s+\eps) -
j(\eps)$ which belongs to $C^1(\R_+,\R)$, and the above computation for the associated
functional $\EE_\eps$ is correct, so that inequality (\ref{2.17}) holds for $\EE$ replaced
by $\EE_\eps$.  Then, writing inequality (\ref{2.17}) for $\EE_\eps$ and fixed $0 \le \phi,\psi \in
L^\infty(\R^N)$, $t
\in (0,1)$ and passing to the limit $\eps \to 0$ we obtain that $\EE$ is convex on $L^\infty(\R^N)$.
Now let us fix $0 \le \phi,\psi \in L^1(\R^N)$, $t \in (0,1)$. 
If $j(\phi)$ or $j(\psi) \notin L^1(\R^N)$ then $t \, \EE(\phi) + (1-t) \, \EE(\psi) = +\infty$ and
the convex inequality (\ref{2.17}) obviously holds. 
In the other case, we have $j(\phi), \, j(\psi) \in L^1(\R^N)$, and we can choose two sequences
$0 \le (\phi_n), \, (\psi_n)$ of $L^\infty(\R^N)$ such that $\phi_n \nearrow \phi$ and $\psi_n
\nearrow \psi$ a.e.. Passing to the limit $\eps \to 0$ in the convex inequality
(\ref{2.17}) written for $\phi_\eps$ and $\psi_\eps$ we get, by the Lebesgue convergence dominated Theorem
and the Fatou Lemma,
\bean
\int_{\R^N} j(t \, \phi + (1-t) \, \psi) &\le &
\mathop{\liminf}_{\eps \to 0} \int_{\R^N} j(t \, \phi_\eps + (1-t) \, \psi_\eps) \\
&\le& t \, \EE(\phi) + (1-t) \, \EE(\psi) + j \Bigl( \int_{\R^N} t \, \phi + (1-t) \, \psi \Bigr),
\eean
which exactly means that $\EE$ is a convex functional in $L^1(\R^N)$. 
Finally, if $0 \le \phi, \, \psi \in L^1(Y \times \R^N)$ and $t \in (0,1)$, then 
$\phi(y,.), \, \psi(y,.) \in L^1(\R^N)$ for almost every $y \in Y$ and, integrating the convex
inequality (\ref{2.17}), we obtain that the functional 
$$
0 \le \phi \in L^1(Y \times \R^N) \mapsto \FF(\phi) = \int_Y \EE(\phi) \, d\nu
$$
is convex. Furthermore, by Fatou Lemma, $\FF$ is l.s.c. for the strong convergence in $L^1$, for the weak $\sigma(L^1,L^\infty)$ convergence and for the  biting
$L^1$-weak convergence, so that (\ref{phiBdd2}) holds.
\qed

\medskip
We introduce a second kind of weak-weak convergence, namely the renormalized convergence, which is the
very natural notion of convergence when we deal with sequences of trace functions, as we will see below. 
We now present  the definition (in a simplified case) and a first elementary result that we
will use in the next subsection. More about the renormalized convergence is presented in section~\ref{sec:WW}. 

\begin{defin}\label{defRenorm}  Let  us define the sequence $(T_M)$ by setting $T_M(s) := s \wedge M = \min(s,M)$ $\forall \, s, \, M \ge 0$. 
We say that a sequence $(\phi_n)$ of $L(Y)$  converges in the renormalized sense (or  r-converges) if there exists a sequence $(\bar T_M)$ of $L^\infty(Y)$ such that
$$
T_M(\phi_n) \rightharpoonup \bar T_M 
\quad \sigma(L^\infty(Y),L^1(Y)) \, \star
\quad\hbox{ and }\quad
\bar T_M \nearrow \phi \quad \hbox{ a.e. in } Y.
$$
\end{defin}

\begin{lem}\label{BitToRenorm} 
For any sequence $(\phi_n)$ of $L(Y)$ and $\phi \in L^0(Y)$ such that
$\displaystyle{\phi_n \bto \phi}$ in the biting $L^1$-weak sense, there exists a subsequence $(\phi_{n'})$ such that 
$\displaystyle{\phi_{n'} \rto \phi}$ in the renormalized sense. 
\end{lem}

\noindent
{\sl Proof of Lemma~\ref{BitToRenorm}.} 
We follow the proof of \cite{8} where that result is established in a $L^1$ framework. 
By assumption, for any $k\in \N$, there exists  a Borel set $A_k$ such that $ \nu \, (Y_k \backslash A_k) < 1/k$ and
$\phi_n \wto \phi$ weakly in $L^1(A_k)$. Thanks to Dunford-Pettis Lemma, there is a function 
$\delta_k:\R_+ \to \R_+$ such that $\delta_k(M) \to 0$ when $M \to + \infty$ and
\beqn\label{BitToR1}
\int_{A_k} \phi_n \, {\bf 1}_{Ê\{Ê\phi_n \ge M \} } \, dy \le \delta_k(M) \quad \forall \, n,M,k \in \N^*.
\eeqn
Moreover, there exists a subsequence $(\phi_{n'})$ of $(\phi_n)$ and a sequence $(\bar T_M)$ of $L^\infty(Y)$ such
that for any $M \in \N$ there holds
$$
T_M(\phi_{n'}) \rightharpoonup \bar T_M  \quad \sigma(L^\infty(Y),L^1(Y)).
$$
We obviously have that $(\bar T_M)$ is an increasing sequence in $L^\infty(Y)$ and $\bar T_M \le \phi$ a.e. because that is true on any $A_k$. 
Observe that 
\beqn\label{BitToR2}
0 \le \phi_n - T_M(\phi_n) \le (\phi_n-M) \, {\bf 1}_{\phi_n \ge M}Ê\quad \hbox{ a.e. in } Y.
\eeqn
Gathering (\ref{BitToR1}) and (\ref{BitToR2}) we get
\bean
\int_{A_k} |Ê\phi - \bar T_M|Ê\, d\nu  
&=& \lim_{n' \to \infty} \int_{A_k} (Ê\phi_{n'} - T_M(\phi_{n'})) \, [\hbox{sign} (Ê\phi - \bar T_M)Ê] Ê\, d\nu  \\
&\le&\liminf_{n' \to \infty} \int_{A_k} \phi_{n'} \, {\bf 1}_{Ê\{Ê\phi_{n'} \ge M \} } \, dy \le \delta_k(M).
\eean
That proves $\bar T_M \to \phi$ a.e. in $Y$ when $M\to\infty$, and then $\displaystyle{\phi_{n'} \rto \phi}$. \qed

%%%%%%%%%%%%%%%%%%%%%%%%%%%%%%%%%%%%%%%%%%%%%%%

\subsection{The trace theorem and the stability result.}

Let us recall the following trace theorem which makes precise the meaning of the trace of a solution. 

\begin{theo}\label{FTeqTraceTheo} {\bf \cite{54}} Let $g\in L^\infty(0,T; L^1(\OO))$ satisfy
$$
\Lambda \, g := \partial_t g + {v} \cdot \nabla_x g = 0 \quad \hbox{in}\quad  \DD'((0,T) \times \OO).
$$
There exists $\gamma g \in L^1_{loc}((0,T) \times \Sigma; d\lambda_2)$ and $g_0 \in L^1(\OO)$ which satisfy the renormalized Green formula
\bear\label{GreenFT}
\int_{0}^T\!\!\int\!\!\!\int_{\OO} \beta(g) \, \Lambda \phi  \, d{v} dx dt 
=  \int_0^T\!\!\int\!\!\!\int_{\Sigma} \beta(\gamma \, g) \, \phi \,\,  n(x) \cdot {v} \,\,
d{v} d\sigma_{\! x} dt - \int\!\!\!\int_{\OO}  \beta( g_0 ) \, \phi \, dx d{v},  
\eear
for all $\beta \in W^{1,\infty}(\R)$ and all test functions $\phi \in \DD([0,T)Ê\times \bar \OO)$, 
as well as for all $\beta \in W^{1,\infty}_{loc}(\R)$, with $\beta' \in L^\infty(\R)$,  and all test functions $\phi \in \DD([0,T)Ê\times \bar \OO)$ such that $\phi = 0$ on $[0,T) \times \Sigma_0$ . 
\end{theo}

\medskip

We may then state our first main result.

\begin{theo}\label{FTeqStab} Let $f_n \in L^\infty(0,\infty;L^1(\OO))$ be a sequence of solutions to the initial boundary value problem (\ref{FTeq})-(\ref{BdaryCond})-(\ref{CIf0}) such that both $(f_n)$ and the trace sequence $(\gamma f_n)$ satisfy the associated  natural a priori bounds: for any $T > 0$ there is a constant $C_T$
\beqn\label{FTeqEstimfn}
\sup_{[0,T]} 
\int \!\! \int_\OO f_n \, \bigl( 1 + |{v}|^2 + |\log \, f_n| \bigr) \, d{v} dx \le C_T
\eeqn
and 
\beqn\label{FTeqEstimgamman}
 \int_0^T \int\!\!\int_\Sigma \gamma f_n \, (1 + |{v}|^2) \, |n(x) \cdot {v}| \, d{v} d\sigma_{\! x} dt
+ \int_0^T \!\! \int_{\partial\Omega} \EE( \gamma_+ f_n ) \, d\sigx dt  \le \alpha^{-1} \, C_T.
\eeqn
On the one hand, there exists $f \in L^\infty(0,\infty; L^1(\OO))$ satisfying (\ref{FTbd1}) and $f_{in} \in L^1(\OO)$ satisfying (\ref{bornef0}) such that, up to the extraction of subsequences, 
\beqn\label{FTeqCvgcefn}
f_n \ \wto \ f \quad \sigma(L^1,L^\infty),  \quad\quad
f_n(0,.) \ \wto \ f_{in} \quad \sigma(L^1,L^\infty) ,
\eeqn
and $f$ is a solution to the free transport equation (\ref{FTeq})-(\ref{CIf0}) with initial datum $f_{in}$.

On the other hand, there exists $\eta_\pm \in L^1(0,T) \times \Sigma_\pm, d\lambda_1)$  for all $T \in (0,\infty)$, which furthermore satisfies (\ref{FTbd1})-(\ref{FTbd2}) (with $\gamma_\pm f$ replaced by $\eta_\pm$),  such that, up to the extraction of a subsequence, for any $T, \eps > 0$ there exists a measurable set $A \subset (0,T) \times \partial\Omega$ such that meas$\, ((0,T) \times \partial\Omega \, \backslash \, A)~<~\eps$  and 
\bear\label{FTeqCvgcegammafn}
&&\gamma_\pm f_n \  \wto \ \eta_\pm \quad\hbox{weakly  in}\quad L^1(A \times \R^N,d\lambda_1), 
\eear
As a consequence, $ \gamma_\pm f = \eta_\pm$ and the reflection boundary condition~(\ref{BdaryCond}) holds. 
\end{theo}

%\begin{rem}\label{FTeqStab} The principle of the proof is the following. 

%- First, from the bounds (\ref{FTeqEstimfn}) and the Dunford-Pettis Lemma, as well as from the  bounds (\ref{FTeqEstimgamman} and the compactness Theorem~\ref{theo1xi} we deduce (\ref{FTeqCvgcefn})  as well as $\gamma_- f_n \wto \eta_-$ in the biting $L^1$-weak sense and $\gamma_+ f_n \wto \eta_+$ in the sense stated in (\ref{FTeqCvgcegammafn}).  

%- Next, Lemma~\ref{BitToRenorm} immediately implies that the above convergences also hold in the renormalized sense. Using the trace Theorem~\ref{FTeqTraceTheo} we then obtain $\gamma_\pm f = \eta_\pm$.   

%- Finally, because $\RR_x$ is continuous in $L^1$ (and thus $\RR_x$ is weakly $L^1$ continuous) we have $\RR \eta_+ = \lim \RR \gamma_+ f_n$ in the sense stated in (\ref{FTeqCvgcegammafn}).  
%\end{rem}

%On the other hand, thanks to (\ref{FTeqCvgcegammafn}), we may clearly pass to the limit $n\to \infty$ in the $L^1$-weak sense in the boundary equation 
%$$
%\gamma_- f_n = \RR ( \gamma_+ f_n) \quad \hbox{on} \quad A \times \R^N,
%$$
%for any set $A$ given by (\ref{FTeqCvgcegammafn}), and we get $\eta_- = \RR ( \eta_+)$. We conclude that the boundary condition (\ref{BdaryCond}) is well fulfilled.
%and the biting Lemma (Theorem~\ref{theobiting} applied to $\psi_n = \gamma_- f_n$) we deduce that $\gamma_- f_n \wto \eta_-$ in the biting $L^1$-weak sense. 

\noindent
{\sl Proof of Theorem~\ref{FTeqStab}. }ÊFirst, from (\ref{FTeqEstimfn}) and the Dunford-Pettis lemma we deduce (\ref{FTeqCvgcefn}). Then, thanks to Lemma~\ref{BitToRenorm}, extracting again a subsequence if necessary, we deduce that $\displaystyle{f_n \rto f}$ or more precisely, there exists two sequences $(\bar T_M)$ and $(\bar T^0_M)$ such that 
\bear\label{stabTL1}
&&T_M(f_{n}) \wto \bar T_M \quad \sigma(L^\infty,L^1) \, \star
\quad\hbox{ and }\quad \bar T_M \nearrow f \quad \hbox{ a.e.}, \\
\label{stabTL2}
&&T_M(f_{n}(0,.)) \wto \bar T_M^0 \quad \sigma(L^\infty,L^1) \, \star
\quad\hbox{ and }\quad \bar T_M^0 \nearrow f_{in} \quad \hbox{ a.e.}.
\eear
\indent
Next,  from (\ref{FTeqEstimgamman}) and Theorem~\ref{theo1xi} (with $\phi_n = \gamma_+ f_n/M$, $j(s) = s \log s$, $\omega(v) = |v|^2$, $d\nu(y) = d\sigma_x dt$, $d\mu_y(v) = |n(x) \cdot v|Ê\, M(v) \, dv$) we deduce that $\gamma_+ f_n \wto \eta_+$ in the sense stated in (\ref{FTeqCvgcegammafn}). That implies that for any $T, \eps > 0$ there exists a measurable set $A \subset (0,T) \times \partial\Omega$ such that meas$\, ((0,T) \times \partial\Omega \, \backslash \, A)~<~\eps$  and
$$
\widetilde{\gamma_+ f_n}\, \wto \, \widetilde{\eta_+}Ê\quad\hbox{weakly  in}\quad L^1(A),
$$
so that $D(\gamma_+ f_n)  \wto D(\eta_+)$ in the sense stated in (\ref{FTeqCvgcegammafn}). 
That also implies that for any $\phi \in L^\infty(\R^N_v)$ 
\bear\nonumber
 \int_{\R^N} \gamma_+ f_n(t,x,v)  \, \phi (R_x \, v) \, n(x) \cdot v \, dv 
\wto \int_{\R^N} \eta_+ (t,x,v)  \, \phi (R_x \, v) \, n(x) \cdot v \, dv  \quad\hbox{weakly  in}\quad L^1(A),
\eear
which means nothing but $L(\gamma_+ f_n)  \wto L(\eta_+)$ in the sense stated in (\ref{FTeqCvgcegammafn}). Gathering these two convergence results, we get   $\gamma_- f_n \wto \eta_-$ in the sense stated in (\ref{FTeqCvgcegammafn}) with $\eta_- := \RR (\eta_+)$. 

Finally,  thanks to Lemma~\ref{BitToRenorm} again, extracting a subsequence if necessary, we deduce that $\displaystyle{\gamma f_n \rto \eta := \eta_+ \, {\bf 1}_{(0,\infty) \times \Sigma_+} + \eta_- \, {\bf 1}_{(0,\infty) \times \Sigma_-} }$ or more precisely, there exists a sequence $(\bar \gamma_M)$ such that 
\beqn\label{stabTL3}
T_M(\gamma f_{n}) \wto \bar \gamma_M \quad \sigma(L^\infty,L^1) \, \star
\quad\hbox{ and }\quad \bar \gamma_M \nearrow \eta \quad \hbox{ a.e.}.
\eeqn
We write then the Green renormalized formula (\ref{GreenFT}) for the free transport equation
$$
\int_{0}^T\!\!\int\!\!\!\int_{\OO}T_M(f_n) \, \Lambda \varphi  \, d{v} dx dt 
=  \int_0^T\!\!\int\!\!\!\int_{\Sigma} T_M(\gamma \, f_n) \, \varphi \,\,  n(x) \cdot {v} \,\,
d{v} d\sigma_{\! x} dt - \int\!\!\!\int_{\OO}  T_M ( f_n(0,.)) \, \varphi \, dx d{v},  
$$
for any $\varphi \in \DD([0,T) \times \bar\OO)$. Using (\ref{stabTL1}), (\ref{stabTL2}) and (\ref{stabTL3}), we may pass twice two the limit in the above equation, first when $n\to \infty$, next when $M\to\infty$, and we get
$$
\int_{0}^T\!\!\int\!\!\!\int_{\OO} f \, \Lambda \varphi  \, d{v} dx dt 
=  \int_0^T\!\!\int\!\!\!\int_{\Sigma} \eta\, \varphi \,\,  n(x) \cdot {v} \,\,
d{v} d\sigma_{\! x} dt - \int\!\!\!\int_{\OO}  f_{in} \, \varphi \, dx d{v}.
$$
In other words, $f$ is a solution to the free transport equation and $\gamma_\pm f = \eta_\pm$ thanks to the trace Theorem~\ref{FTeqTraceTheo}. 
We conclude by gathering that information with the equation satisfied by $\eta_\pm$. \qed

%%%%%%%%%%%%%%%%%%%%%%%%%%%%%%%%%%%%%%%%%%%%%%%%
\section{On the convergence in the renormalized sense.}\label{sec:WW}
\setcounter{equation}{0}
\setcounter{theo}{0}
%%%%%%%%%%%%%%%%%%%%%%%%%%%%%%%%%%%%%%%%%%%%%%%%

\subsection{Basic properties.}

We present the main basic properties concerning the notion of convergence in the renormalized sense. 
More about renormalized convergence is set out in the appendix section. In that section the framework and notations are the same as those of subsection 2.2, and again, we only deal with non-negative functions of $L = L(Y)$, but we do not specify it anymore.

\begin{defin}\label{defWW1} We say that $\alpha$ is a renormalizing function if $\alpha \in C_b(\R)$
is increasing and  $0 \le \alpha(s) \le s$ for any $s \ge 0$. 
We say that $(\alpha_M)$ is a renormalizing sequence if $\alpha_M$ is a renormalizing function for
any $M \in \N$ and $\alpha_M(s) \nearrow s$ for all $s \ge 0$ when $M \nearrow \infty$.
Given any renormalizing sequence $(\alpha_M)$, we say that $(\phi_n)$ $(\alpha_M)$-renormalized
converges to $\phi$ (or we just say that $(\phi_n)$ r-converges to $\phi$) if there exists a sequence $(\bar\alpha_M)$ of $L^\infty(Y)$ such that
$$
\alpha_M(\phi_n) \rightharpoonup \bar \alpha_M 
\quad \sigma(L^\infty(Y),L^1(Y)) \, \star
\quad\hbox{ and }\quad
\bar\alpha_M \nearrow \phi \quad \hbox{ a.e. in } Y.
$$
Notice that the renormalized convergence as defined in definition~\ref{defRenorm} is nothing but the $(T_M)$-renormalized convergence. \end{defin}

\begin{prop}\label{propWW1}
1. The $(\alpha_M)$-renormalized limit in the definition~\ref{defWW1} does not depend on the renormalizing
sequence $(\alpha_M)$, but only on the sequence $(\phi_n)$. 
In other words, given two renormalizing sequences $(\alpha_M)$ and $(\beta_M)$, if $(\phi_n)$ $(\alpha_M)$-renormalized converges to $\phi^\alpha$ and $(\beta_M)$-renormalized converges to $\phi^\beta$ then $\phi^\alpha = \phi^\beta$. 

2. For any sequence $(\phi_n)$ of $L$ there exists a subsequence $(\phi_{n'})$ of $(\phi_n)$ and a function $\phi \in L$ such that 
$(\phi_{n'})$ $(\alpha_M)$-renormalized converges to $\phi$ for any renormalizing sequence $(\alpha_M)$.

3. A sequence $(\phi_n)$ which converges to $\phi$ a.e. or strongly  in $L^p$, $p \in [1,\infty]$, also r-converges to $\phi$. 
From a sequence $(\phi_n)$ which converges to $\phi$ weakly  in $L^p$, $p \in [1,\infty]$, or in the biting $L^1$-weak sense, we may extract a subsequence  $(\phi_{n'})$ which r-converges to $\phi$.
\end{prop}

\begin{rem}\label{remWW1}
1.  The definition of the $(\alpha_M)$-renormalized convergence with $\alpha_M \not=
T_M$ is important in order to obtain the renormalized convergence of the trace functions sequence in Theorem~\ref{th:rcvgeVFP}. 
Indeed, $T_M$ is not smooth enough in order to be taken as a renormalizing function for the
VFP equation and we have to introduce the ``smooth" renormalizing functions $\alpha :=
\Phi_{M,\theta}$.

2.  Because of Proposition~\ref{propWW1}  we will often make the abuse of language by not specifying the renormalizing sequence $(\alpha_M)$ used to define the $(\alpha_M)$-renormalized convergence and by saying that $(\phi_n)$ r-converges (to $\phi$) when it is only a subsequence of $(\phi_n)$ which r-converges (to $\phi$).

3.  Let us notice that in general we can not exclude that the limit $\phi \equiv + \infty$, since for instance the sequence $(\phi_n)$ defined by $\phi_n =  n$ belongs to $L$ and r-converges to $\phi \equiv \infty$. 

\end{rem}

\noindent
{\sl Proof of the Proposition~\ref{propWW1}.}
{\sl Step 0. }ÊWe first claim that for any sequence $(\phi_n)$ of $L$ and any renormalizing sequence $(\alpha_M)$ there exists a
subsequence $(\phi_{n'})$ of $(\phi_n)$ and $\phi \in L$ such that $(\phi_{n'})$ $(\alpha_M)$-renormalized converges to $\phi$. 
Indeed, for any $M$ we can find a subsequence $(n^M_k)_k$ and $\bar\alpha_M \in L^\infty$ such that $\alpha_M(\phi_{n^M_k}) \wto \bar\alpha_M$ weakly in $L^\infty$. By a diagonal process we can obtain a unique subsequence $(n')$ such that the above weak convergence holds for any $M \in \N$. Furthermore, since $(\alpha_M)$ is increasing, we get that $(\bar\alpha_M)$ is an increasing
sequence of non-negative measurable functions, so that it converges to a limit $\phi \in L$.

\smallskip\noindent
{\sl Step 1.} Assume that for a renormalizing sequence $(\alpha_K)$ we have
$\displaystyle{ \alpha_K(\phi_n) \wto \bar\alpha_K \nearrow \psi}$.
Thanks to step 0, there exits a sub-sequence $(\phi_{n'})$, a sequence $\bar T_M \in L^\infty$
and a function $\phi \in L$ such that $\displaystyle{T_M(\phi_{n'}) \wto \bar T_M \nearrow \phi}$.
It is clear that $\forall K,M \in \N$ $\forall \eps > 0$ there is $k_{M,\eps}, m_K \in\N$ such that
$\alpha_K \le T_{m_K}$ and $T_M \le \alpha_{k_{M,\eps}} + \eps$.
Therefore, writing that $\alpha_K(\phi_n) \le T_{m_K}(\phi_n)$ and 
$T_M(\phi_n) \le \alpha_{k_{M,\eps}}(\phi_n) + \eps$, and passing to the limit $n \to + \infty$, we
get
$$
\bar \alpha_K \le \bar T_{m_K} \le \phi \quad\hbox{ and }\quad
\bar T_M \le \bar\alpha_{k_{M,\eps}} + \eps \le \psi +Ê\eps.
$$
Then passing to the limit $M,K \nearrow\infty$ we obtain that $\psi \le \phi \le \psi + \eps$ for any $\eps > 0$,  and finally passing to the limit $\eps \to 0$ we conclude that $\psi = \phi$.

\smallskip\noindent
{\sl Step 2.}
Let us remark that the class of renormalizing functions is separable for the uniform norm of $C(\R_+)$.
For instance, the family $\AA = \{ \alpha^k \}$ of functions $\alpha$ such that 
$$
0 \le \alpha (s) \le s \quad\hbox{ and }\quad
\alpha'(s) = \sum_{j = 1}^J \theta_j \, {\bf 1}_{[a_j,a_{j+1}[}(s), \quad a_j, \, \theta_j \in \Q_+
$$
is countable and dense. By a diagonal process and thanks to step 0, we can find a subsequence $(\phi_{n'})$ in such a
way that for any $\alpha \in \AA$ there exists $\bar\alphaÊ\in L^\infty$ such that
$\alpha(\phi_{n'}) \wto \bar\alpha$. Let us fix now $(\beta_M)$ a renormalizing sequence. 
On one hand, for any $M$ there exists a sequence $(\alpha_k)$ of $\AA$ such that $\alpha_k  \le \beta_M \le \alpha_k + 1/k$ for any $k \in \N$ and $\alpha_k \nearrow \beta_M$. We already know that $\alpha_k(\phi_{n'}) \wto
\bar\alpha_k$. Since $(\bar\alpha_k)$ is not decreasing, it converges a.e., and we set 
$\beta_M^* = \lim \bar\alpha_k$.
On the other hand, thanks to Step 0, there exists a subsequence $(\phi_{n''})$ and a function
$\bar\beta_M$ such that $\beta_M(\phi_{n''}) \wto \bar\beta_M$. That implies 
$\bar\alpha_k \le \bar\beta_M \le \bar\alpha_k + 1/k$. Passing to the limit $k \to \infty$,
 we get $\bar\beta_M =\beta_M^*$. Therefore, by uniqueness of the limit, it is the all sequence
$\beta_M(\phi_{n'})$ which converges to $\beta^*_M$. 
Finally, thanks to the usual monotony argument we deduce that $\phi_{n'}$ converges in the $(\beta_M)$-renormalized sense and its limit is necessary $\phi$ thanks to Step 1. 

\smallskip\noindent
{\sl Step 3.}  If $\phi_n \to \phi$ a.e. then clearly $\alpha_M(\phi_n) \wto \alpha_M(\phi)$ $L^\infty$-weak and
$\alpha_M(\phi) \nearrow \phi$ for any renormalizing sequence $(\alpha_M)$, so that $\displaystyle{\phi_n \rto \phi}$. If $(\phi_n)$ converges strongly or weakly in $L^p$, $p\in [1,\infty]$, then it obviously converges in the biting $L^1$-weak sense and we may apply Lemma~\ref{BitToRenorm}. 
 \qed

\medskip
Let us now define the limit superior and the limit inferior in the renormalized sense.

\begin{defin}\label{defWW2}
Let $(\phi_n)$ be a sequence of $L$. Consider $I$ the set of all
the increasing applications $\imath : \N \to \N$ such that the subsequence
$(\phi_{\imath(k)})_{k\ge 0}$ of $(\phi_n)_{n \ge 0}$ converges in the renormalized sense and note 
$\phi_\imath = \rlim \phi_{\imath(k)}$. Thanks to the Proposition~\ref{propWW1}.2, we know that $I$ is
not empty. We defined the limit superior and the limit  inferior of
$(\phi_n)$ in the renormalized sense by
$$
\rlimsup \phi_n := \sup_{ \imath \in I} \phi_\imath
\quad\hbox{ and }\quad
\rliminf \phi_n := \inf_{ \imath \in I} \phi_\imath.
$$
It is clear that if $\rlimsup \phi_n =\rliminf \phi_n$ then $(\phi_n)$ r-converges (up to the extraction of 
a subsequence).
\end{defin}

\begin{prop}\label{propWW3}
1. If $\displaystyle{\phi_n \rto \phi}$, $\displaystyle{\psi_n \rto \psi}$ and $\lambda_n \to
\lambda$ in $\R_+^\star$ then $\displaystyle{\phi_n + \lambda \, \psi_n \rto \phi + \lambda \,
\psi}$.

2. Let $\displaystyle{\phi_n \rto \phi}$ and $\beta$ be a non-negative and concave function then
 $\beta(\phi) \ge \rlimsup \beta(\phi_n)$.

3. Let $\beta$ be a strictly concave function, and $(\phi_n)$ be a sequence such that
$\displaystyle{\phi_n \rto \phi}$ and $\beta(\phi) \le \rliminf
\beta(\phi_n)$ then, up to the extraction a subsequence, $\phi_n \to \phi$ a.e. in $Y$.

4. Let $\displaystyle{\phi_n \rto \phi}$ and $S$ be a bounded and non-negative operator of
$L^1$ then $S \, \phi \le \rliminf S \, \phi_n$.
\end{prop}

\medskip\noindent
{\sl Proof of the Proposition~\ref{propWW3}.} 
{\sl Step 1.}  From the elementary inequality
$$
\forall \, a, b, M \ge 0 \qquad M \wedge (a+b) \le M \wedge a + M \wedge b \le (2 \, M) \wedge (a+b),
$$
we deduce
$$
\wlim [ M \wedge (\phi_n+\psi_n)] \le \wlim [M \wedge \phi_n] + \wlim [M \wedge \psi_n]
\le \wlim [ (2 \, M) \wedge (\phi_n+\psi_n)]
$$
so that $\rlim(\phi_n + \psi_n) = \phi + \psi$. Next, from the elementary identity
$$
\forall \, a, b, M \ge 0 \qquad 
(a \, b) \wedge M = a \, \bigl( b \wedge (M/a) \bigr)
$$
and because for any $\eps > 0$ there holds $0 < \lambda-\eps \le \lambda_n \le \lambda+\eps$  for $n$ large enough, we have
$$
(\lambda-\eps) \, \Bigl[\phi_n \wedge {M \over \lambda-\eps} \Bigr]
\le (\lambda_n \, \phi_n) \wedge M
\le (\lambda+\eps) \, \Bigl[\phi_n \wedge {M \over \lambda+\eps} \Bigr].
$$
We deduce that for a subsequence $(\lambda_{n'} \, \phi_{n'})$ 
$$
\forall \, \eps > 0 \qquad (\lambda-\eps) \, \bar T_{M \over \lambda-\eps} 
\le \mathop{\wlim}_{n' \to \infty} T_M (\lambda_{n'} \, \phi_{n'})
\le (\lambda+\eps) \, \, \bar T_{M \over \lambda+\eps} ,
$$
so that, passing to the limit $\eps \to 0$ and using that $T_{M/(\lambda+\eps)} \le T_{M/(\lambda-\eps)}$,
$$
\lambda \, \bar T_{M} =  \mathop{\wlim}_{n' \to \infty} T_M (\lambda_{n'} \, \phi_{n'}).
$$
Passing to the limit $M \to \infty$, we conclude that $\rlim (\lambda_n \, \phi_n)=  \lambda \, \phi$. 

\smallskip\noindent{\sl Step 2.} 
We know that
$$
\beta(s) = \inf_{\ell \ge \beta} \ell(s),
$$
where the inf is taken over all real values affine functions $\ell(t) = a \, t + b$  which satisfy $a,b \ge
0$ and $\beta(t) \le \ell(t)$ for any $t \ge 0$. 
Furthermore, for any $\ell$ and $M$, there clearly exists $K_M$ such that
$$
T_M(\ell(s)) \le \ell(T_K(s))
\qquad\hbox{and}\qquad
\ell(T_M(s)) \le T_K(\ell(s))
\qquad\hbox{for all}\qquad
K \ge K_M, \ s \ge 0.
$$
We deduce that for any $\ell \ge \beta$, we have
$$
T_M(\beta(\phi_n)) \le \ell(T_K(\phi_n)).
$$
Therefore, we get 
$$
\limsup_n T_M(\beta(\phi_n)) \le \ell(\lim_n T_K(\phi_n)) \le \ell (\phi)
$$
and finally
$$
\limsup_n T_M(\beta(\phi_n)) \le \beta (\phi) \qquad\hbox{for any } M,
$$
which exactly means that $\rlimsup \beta(\phi_n) \le \beta (\phi)$.

\smallskip\noindent{\sl Step 3.} 
For any subsequence $(n')$ such that $\beta(\phi_{n'})$, $\beta(\phi_{n'}/2 + \phi/2))$ and 
$\beta(\phi_{n'}/2 + \phi/2)) - \beta(\phi_{n'})/2 - \beta(\phi)/2 \ge 0$ converge in the
renormalized sense, we have 
$$
\rlim \bigl[ \beta \bigl( {\phi_{n'} + \phi \over 2} \bigr) - {\beta(\phi_{n'}) \over 2}
- {\beta(\phi) \over 2} \bigr] + {\beta(\phi) \over 2} + \rlim  {\beta(\phi_{n'}) \over 2} 
= \rlim \beta \bigl( {\phi_{n'} + \phi \over 2} \bigr),
$$
thanks to step 1. As a consequence,  we get
\bean
0 &\le& \rlim \bigl[ \beta \bigl( {\phi_{n'} + \phi \over 2} \bigr) - {\beta(\phi_{n'}) \over 2}
- {\beta(\phi) \over 2} \bigr] \\
&=& \rlim \beta \bigl( {\phi_{n'} + \phi \over 2} \bigr)  - {\beta(\phi) \over 2} - \rlim 
{\beta(\phi_{n'}) \over 2}  \\
&\le& \beta(\phi) - {\beta(\phi) \over 2} - {\beta(\phi) \over 2} = 0,
\eean
thanks to step 2 and because  $\displaystyle{\phi_{n'}/2 + \phi/2 \rto \phi}$. Therefore, for any
$k$, we have
\bean
0 &\le& \lim_{n \to \infty} \int_{Y_k} T_1 \Bigl( 
\beta \bigl( {\phi_{n'} + \phi \over 2} \bigr) - {\beta(\phi_{n'}) \over 2} - {\beta(\phi) \over 2}
\Bigr) \, d\nu  \\
&\le& \int_{Y_k}
\rlimsup \bigl[ \beta \bigl( {\phi_{n'} + \phi \over 2} \bigr) - {\beta(\phi_{n'}) \over 2}
- {\beta(\phi) \over 2} \bigr] \, d\nu = 0,
\eean
so that, up to extraction a subsequence,
$$
\beta \bigl( {\phi_{n'} + \phi \over 2} \bigr) - {\beta(\phi_{n'}) \over 2} - {\beta(\phi) \over 2} \to 0
\quad \hbox{a.e. on } Y \quad\hbox{ and }\quad
\phi_{n'} \to \phi \quad \hbox{a.e. on } Y. 
$$

\smallskip\noindent{\sl Step 4.} 
Fix $\chi \in C_c(Y)$, the space of continuous functions on $Y$ with compact support, such that $0 \le \chi \le 1$.
Since $T_M(\phi_n) \, \chi \rightharpoonup \bar T_M \, \chi $ weakly in $L^1$, we have 
\beqn\label{6.4}
S( T_M(\phi_n) \, \chi) \rightharpoonup S (\bar T_M \, \chi)
\quad\hbox{ weakly in } L^1.
\eeqn
We deduce, using $T_K ( S ( T_M(\phi_n) \, \chi) ) \le T_K (S ( \phi_n))$ and Proposition~\ref{propWW1}.3 that
$$
S (\bar T_M \, \chi) = \mathop{\rliminf}_{n \to \infty} S ( T_M(\phi_n) \, \chi) \le 
\mathop{\rliminf}_{n \to \infty} S ( \phi_n).
$$ 
We conclude letting $\chi \nearrow 1$ and $M \to +\infty$. \qed

%%%%%%%%%%%%%%%%%%%%%%%%%%%%%%%%%%%%%%%%%%%%%%%%

\subsection{From renormalized convergence to weak convergence.}

We give now a kind of extension of the biting Lemma in the $L^0$ framework.

\begin{defin}\label{defBit2}
We say that a sequence $(\psi_n)$ is asymptotically bounded in $L^0(Y)$ if  for any $k \in \N$ there
exists $\delta_k:\R_+ \to \R_+$ such that $\delta_k(M) \searrow 0$ when $M \nearrow +\infty$ and
for any $M$ there is $n_{k,M}$ such that
\beqn\label{asymBdL0}
\hbox{meas} \, \{ y \in Y_k, \ \psi_n(y) \ge M \}Ê\le \delta_k(M) \qquad \forall  k \in \N, \
\forall n \ge n_{k,M}.
\eeqn
\end{defin}

\begin{theo}\label{theoBit2}
Let $(\psi_n)$ be a sequence of $L^0(Y)$ which r-converges to $\psi$ with $\psi \in L^0(Y)$. 
Then  $(\psi_n)$ is asymptotically bounded in $L^0(Y)$ and there exists a subsequence $(\psi_{n'})$ which b-converges to $\psi$. \end{theo}

\begin{rem}\label{remBit2}
In the $L^1$ framework, J. Ball \& F. Murat \cite{8} have already proved that the biting $L^1$-weak
convergence implies, up to the extraction of a subsequence, the convergence in the renormalized
sense, as it has been recalled and extended to the $L^0$ framework in Lemma~\ref{BitToRenorm}. 
As a consequence, combining Ball \& Murat's result with Theorem~\ref{theoBit2}, 
we get the equivalence between the biting $L^1$-weak convergence
and the renormalized convergence. More precisely, considering a sequence $(\psi_n)$ of $L(Y)$, it is
equivalent to say that, up to the extraction of a subsequence, 
\bear\label{2.7}
&&\psi_n \bto \psi \ \hbox{ in the biting } L^1 \hbox{-weak sense (so that }
\psi \in L^0(Y) \hbox{)}, \\ \label{2.8}
&&\psi_n \rto \psi \ \hbox{ in the renormalized sense and } \psi \in L^0(Y).
\eear
Furthermore, in both cases, the full sequence $(\psi_n)$ is asymptotically  bounded in $L^0$. 
Again, we refer to the appendix where some complements about r-convergence and b-convergence are given. 
\end{rem}

\noindent
{\sl Proof of Theorem~\ref{theoBit2}.} {\sl Step 1. Proof of the asymptotic boundedness in $L^0$}. We argue by
contradiction. For an arbitrary $\eps > 0$ we know that there exists $B \subset Y_k$ such that
$\nu(Y_k \backslash B) < \eps/2$ and $\psi \in L^1(B)$. If there is no $m \in \N$ such that meas$\,
\{y \in B, \,  \psi_n(y) \ge m \} < \eps/2$ for all $n$ large enough, this means that there exists an increasing
sequence $(n_m)$ such that
$$
\hbox{meas} \, \{ y \in B, \, \psi_{n_m}(y)  \ge m \} \ge \eps/2 \qquad \forall m \ge 0.
$$
Therefore, for any $\ell \in \N$ and any $m \ge \ell$ we have 
$$
\int_B T_\ell(\psi_{n_m}) \ge \ell \, \hbox{meas} \, \{ y \in B, \, \psi_{n_m}(y) \ge  \ell \}  
\ge \ell \,  \hbox{meas} \, \{ y \in B, \, \psi_{n_m}(y) \ge m \}  \ge  \ell \, {\eps \over 2},
$$
and passing to the limit $m \to \infty$, we get
$$
\int_B \psi \ge \int_B \mathop{\wlim}_{m \to \infty}Ê T_\ell(\psi_{n_m})  \ge \ell \,  {\eps \over 2}
\qquad \forall \ell \ge 0.
$$
Letting $\ell \nearrow \infty$ we get a contradiction with the fact that $\psi \in L^1(B)$. As a conclusion, we have proved that 
for any $\eps > 0$ there exists $m_\eps$ and $n_\eps$ such that meas$\, \{ y \in Y_k, \, \psi_n \ge m_\eps \} < \eps$ for any $n \ge n_\eps$, and (\ref{asymBdL0}) easily follows. 

\smallskip\noindent 
{\sl Step 2. Proof of the convergence in the biting $L^1$-weak sense.}  As in Step 1, for any $k \in \N$ we can
choose $B$ such that $\nu(Y_k \backslash B) < 1/3 k$ and $\psi \in L^1(B)$. Setting 
$\displaystyle{\int_{B} \psi \, dy = C_0}$, we construct a sequence $(n_\ell)$ such
that 
\beqn\label{TlB}
\int_{B} T_\ell(\psi_{n_\ell}) \, dy \le C_0 + {1 \over \ell}.
\eeqn
From (\ref{TlB}), Theorem~\ref{theobiting} (biting Lemma) and Lemma~\ref{BitToRenorm}, we may extract
a subsequence, still denoted by $(\psi_{n_\ell})$, which b-converges and r-converges to a limit denoted by $\psi^* \in L^1(B)$. 
On the one hand, for any $M \in \N$ we have $T_M (\psi_{n_\ell}) \le T_\ell(\psi_{n_\ell})$ for $\ell \ge M$ so
that, passing to the limit $\ell \to \infty$, we get $\wlim T_M (\psi_{n_\ell}) \le \psi^*$ and thus $\psi \le \psi^*$. On the other hand, 
from Theorem~\ref{theobiting} (biting Lemma) again, we have 
$\|Ê\psi^* \|_{L^1} \le \liminf \|ÊT_\ell(\psi_{n_\ell}) \|_{L^1} \le C_0 = \|Ê\psi \|_{L^1}$. 
Gathering these two inequalities, we have proved 
$$
T_\ell(\psi_{n_\ell}) \wto \psi \quad\hbox{ weakly in }ÊL^1(B).
$$

\smallskip\noindent
Furthermore, since $(\psi_n)$ is asymptotically bounded in $L^0(Y)$ we have, up to the extraction of
a subsequence again,
$$
\hbox{meas} \{ \psi_{n_\ell} \not= T_\ell(\psi_{n_\ell}) \} = \hbox{meas} \{ \psi_{n_\ell} > \ell
\} 
\le \delta_k(\ell) \mathop{\longrightarrow}_{\ell \to \infty} 0.
$$
Therefore, we can choose an other subsequence, still noted $(\psi_{n_\ell})$, such that
$Z_L := \{ \forall \ell \ge L \ / \ \psi_{n_\ell} \not= T_\ell(\psi_{n_\ell}) \}$ satisfies
$$
\hbox{meas} (Z_L) \le \sum_{\ell \ge L} \hbox{meas} \{ \psi_{n_\ell} > \ell \} 
\mathop{\longrightarrow}_{L \to \infty} 0.
$$
Finally, choosing $L$ large enough such that meas$\,(Z_L) < 1/ 3 k$ and setting 
$A_k := B \cap Z_L^c$, we have $|Y_k \backslash A| < 1/k$, $\psi_{n_\ell} \in L^1(A)$ for all 
$\ell \ge L$ and
$$
\psi_{n_\ell} = T_\ell (\psi_{n_\ell}) \wto \psi \quad\hbox{ weakly in }ÊL^1(A).
$$
We conclude thanks to a diagonal process. \qed

\medskip\noindent
A simple but fundamental consequence of Theorem~\ref{theo1xi} and Theorem~\ref{theoBit2} is the following. 

\begin{theo}\label{theoxi2} Consider a function $m : \R^N \to \R$ and a family of measures
$d\varpi_y$ on $\R^N$ such that
$$
\int_{\R^N} m({v}) \, d\varpi_y({v})  = 1, \quad
\int_{\R^N} m({v})^{1/4} \, d\varpi_y({v}) \le C_4 \ \forall y
\quad\hbox{and}\quad m(0) \ge m({v}) \mathop{\longrightarrow}_{|{v}| \to \infty} 0.
$$
Let $(\phi_n)$ be a sequence of $L^0(Y \times \R^N)$ which satisfies
$$
\int_Y  \EE \Bigl( {\phi_n(y,.) \over m(.) }Ê\Bigr) \, d\nu(y) \le C_1 < \infty,
$$
with $\EE$ just like in Theorem~\ref{theo1xi} with $d\mu_y({v}) = m({v}) \, d\varpi_y({v})$,
and assume that 
\beqn\label{psinrtopsi}
\psi_n(y) := \int_{\R^N} \phi_n(y,{v}) \, d\varpi_y({v}) \ \rto \ \psi 
\quad \hbox{with}\quad \psi \in L^0(Y).
\eeqn
Then, there exists $\phi \in L^1(Y \times \R^N,d\nu d\varpi)$ and a subsequence $(\phi_{n'})$
such that for every $k \in \N$ we can find $A_k \subset Y_k$ in such a way
that $(A_k)$ is increasing, $\nu(Y_k \backslash A_k) < 1/k$ and
$$
\phi_{n'} \rightharpoonup \phi \quad\hbox{ weakly in }\quad L^1(A_k \times \R^N,d\nu d\varpi).
$$
As a consequence $\displaystyle{ \psi = \int_{\R^N} \phi \, d\varpi }$ and 
$\EE ( \phi/ m) \in L^1(Y)$. 
\end{theo}

\medskip\noindent{\sl Proof of Theorem~\ref{theoxi2}.} 
From (\ref{psinrtopsi}), Theorem~\ref{theoBit2} and Definition~\ref{defbiting} we know that there exists a subsequence $(\psi_{n'})$ such that for every 
$k \in \N$ we can find $A = A_k \subset Y_k$ satisfying $(A_k)$ is increasing, 
$\nu(Y_k \backslash A_k) < 1/k$ and
$$
\psi_{n'} \quad\hbox{ is weakly compact in }\quad L^1(A).
$$
Next, we come back to estimate (\ref{Psi+}) in the proof of Theorem~\ref{theo1xi}, which written with the new
notation, becomes
\beqn\label{2.19}
\int \!\!\! \int_{A \times \R^N} \phi_{n'} \, \Xi \bigl( {\phi_{n'} \over m({v})} \bigr) \,
d\varpi_y d\nu \le C_3,
\eeqn
where we have set $ \Xi(s) := \Psi^+(s) / s$. Of course, we can assume without loss of generality
that $\Xi$ is not decreasing, $\Xi(s) \nearrow \infty$ when $s \nearrow \infty$ and
$\Xi(s) \le s^{1/2}$. From (\ref{2.19}) we deduce
\beqn\label{2.20}
\int \!\!\! \int_{A \times \R^N} \phi_{n'} \, \Xi \bigl( {\phi_{n'} \over m(0)} \bigr) \, d\varpi_y
\, d\nu \le C_3,
\eeqn
as well as 
\bear\label{2.21}
&&\int \!\!\! \int_{A \times \R^N}  \!\!\! \!\!\!
 \phi_{n'} \, \Xi ( m({v})^{-1/2} ) \, d\varpi_y \, d\nu \le \\ \nonumber
&&\qquad
\le \int \!\!\! \int_{A \times \R^N}  \!\!\! \!\!\! \phi_{n'} \, \Xi ( m({v})^{-{1 / 2}} ) \,
\bigl(  {\bf 1}_{ \{ \phi_{n'} \le m({v})^{1/2} \} } + {\bf 1}_{ \{ \phi_{n'} \ge m({v})^{1/2} \} }
\bigr) \, d\varpi_y \, d\nu \\ \nonumber
&&\qquad
\le \int \!\!\! \int_{A \times \R^N}  \!\!\! \!\!\! m({v})^{1/4}
\, d\varpi_y \, d\nu +
\int \!\!\! \int_{A \times \R^N}  \!\!\! \!\!\! \phi_{n'} \, \Xi \bigl( {\phi_{n'} \over m({v})}
\bigr) \, d\varpi_y \, d\nu \le C_4 \, |Y_k| + C_3.
\eear
Gathering (\ref{2.20}) and (\ref{2.21}), we deduce thanks to the Dunford-Pettis Lemma that $(\phi_{n'})$
belongs to a weak compact set of $L^1(A \times \R^N,d\nu d\varpi)$, and we conclude as in the end
of the proof of Theorem~\ref{theo1xi}. \qed

%%%%%%%%%%%%%%%%%%%%%%%%%%%%%%%%%%%%%%%%%%%%%%%%
\section{Trace theorems for solutions of the Vlasov-Fokker-Planck equation.}\label{sec:Trace}
\setcounter{equation}{0}
\setcounter{theo}{0}
%%%%%%%%%%%%%%%%%%%%%%%%%%%%%%%%%%%%%%%%%%%%%%%%

%%%%%%%%%%%%%%%%%%%%%%%%%%%%%%%%%%%%%%%%%%%%%%%

\subsection{Statement of the trace theorems}

In this section we recall the trace results established in \cite{53}, \cite{54} for the
Vlasov equation (which corresponds to the case $\nu=0$ in the Theorem below) 
and we extend them to the VFP equation. Given a vector field $E=E(t,x,{v})$, a source term $G = G(t,x,{v})$, a constant
$\nu \ge 0$ and a solution $g = g(t,x,{v})$ to the Vlasov-Fokker-Planck equation 
\beqn\label{3.1}
\Lambda_E \, g = {\partial g \over \partial t}  +  {v} \cdot \nabla_x  g + E \cdot \nabla_{v} g 
- \nu \, \Delta_{v} g= G \ \hbox{ in } (0,T) \times \OO,
\eeqn
we show that $g$ has a trace $\gamma g$ on the boundary $(0,T) \times \Sigma$ and a trace
$\gamma_t g$ on the section $\{ t \}Ê\times \OO$ for any $t \in [0,T]$. These trace functions
are defined thanks to a Green renormalized formula. We write indifferently $\gamma_t g = g(t,.)$. 

\smallskip
The meaning of equation (\ref{3.1}) is of two kinds.
In the first case, we assume that $g \in L^\infty(0,T;$ $L^p_{loc}(\bar \OO))$, with $p \in
[1,\infty]$, is a solution of (\ref{3.1}) in the sense of distributions, i.e.,
\beqn\label{3.2}
\int_0^T\!\!\!\int\!\!\!\int_\OO (g \, \Lambda^\star_E \phi + G \, \phi) \, d{v} dx dt = 0, 
\eeqn
for all test functions $\phi \in \DD((0,T) \times \OO)$, where we have set
$$
\Lambda^\star_E \, \phi = {\partial \phi \over \partial t}  +  {v} \cdot \nabla_x  \phi + E \cdot
\nabla_{v} \phi  + \nu \, \Delta_{v} \phi + (\hbox{div}_{v} E) \, \phi.
$$ 
In this case we assume
\beqn\label{(3.4)}
E \in  \,\, L^1 \bigl( 0,T;W^{1,p'}_{loc}(\bar\OO) \bigr), \quad
\hbox{div}_{v} E \in L^1 \bigl( 0,T;L^{p'}_{loc}(\bar\OO) \bigr), \quad
G \in L^1_{loc}([0,T] \times \bar \OO), 
\eeqn
where $p' \in [1,\infty]$ stands for the conjugate exponent of $p$, given by $1/p +1/p' = 1$,
and we make one of the two additional hypothesis
\beqn\label{(3.5)}
\nu \int_0^T \!\! \int_{\OO_R} | \nabla_{v} g |^2 \, d{v} dx dt \le C_{T,R} \\
\eeqn
or
\beqn\label{(3.6)}
\nu \int_0^T \!\! \int_{\OO_R} | \nabla_{v} g |^2 \, 
{\bf 1}_{ \{ M \le |g| \le M+1 \} } \, d{v} dx dt \le C_{T,R} \quad \forall M \ge 0. 
\eeqn

\begin{rem}\label{rem3.1}
The bound (\ref{(3.6)}) is the natural bound that appears when we consider, for example, the initial
value problem with initial datum $g_0 \in L^p(\OO)$ when $\Omega = \R^N$ or when $\Omega$ is an
open subset of $\R^N$ and specular reflections are imposed at the boundary.
\end{rem}

\smallskip
In the second case, we assume that $g$ is a renormalized solution of (\ref{3.1}). In order to make precise
the meaning of such a solution, we must introduce some notations. We denote by $\BB_1$ the class
of functions $\beta \in W^{2,\infty}(\R)$ such that
$\beta'$ has a compact support and by $\BB_2$ the class of functions $\beta \in
W^{2,\infty}_{loc}(\R)$ such that $\beta''$ has a compact support.
Remark that for every $u \in L(Y)$ and $\beta \in \BB_1$ one has $\beta(u) \in L^\infty(Y)$. 
We shall write $g \in C([0,T];L(\OO))$ if $\beta(g) \in C([0,T];L^1_{loc}(\bar\OO))$ for every
$\beta \in \BB_1$.

We say that $g \in L((0,T) \times \OO)$ is a renormalized solution of (\ref{3.1}) if for all
$\beta \in \BB_1$ we have
\beqn\label{3.7}
E \in L^1 \bigl( 0,T; W^{1,1}_{loc}(\bar\OO) \bigr), 
\ \beta'(g) \, G \in L^1_{loc}([0,T] \times \bar \OO), 
\ \nu \,  \beta''(g) \, |\nabla_{v} \, g |^2 \in L^1_{loc}([0,T] \times \bar \OO),
\eeqn
and $\beta(g)$ is solution of
\beqn\label{3.8}
\Lambda_E \, \beta (g) = \beta'(g) \, G - \nu \, \beta''(g) \, |\nabla_{v} \, g |^2 
\ \hbox{ in } \DD'((0,T) \times \OO). 
\eeqn

\smallskip

We can now state the trace theorems for the Vlasov-Fokker-Planck equation (\ref{3.1}). 

\begin{theo}\label{theoTraceLinfty}
{\bf (The case $p =\infty$)}. 
Let $g \in L^\infty([0,T]Ê\times \OO)$ be a solution of equation (\ref{3.2})-(\ref{(3.4)})-(\ref{(3.5)}). 
There exists $\gamma g$ defined on $(0,T) \times \Sigma$ and for every $t \in [0,T]$ there exists $\gamma_t g \in L^\infty(\OO)$ such that
\beqn\label{3.9}
\gamma_t g \in C([0,T];L^a_{loc}(\bar\OO)) \quad \forall a \in [1,\infty)
\quad \hbox{and } \quad
\gamma \, g \in L^\infty((0,T)Ê\times \Sigma),
\eeqn
and the following Green renormalized formula
\bear\label{3.10}
&&\int_{t_0}^{t_1}\!\!\int\!\!\!\int_{\OO} \bigl( \beta(g) \, \Lambda^\star_E \phi 
+ (\beta'(g) \, G - \nu \, \beta''(g) \, |\nabla_{v} \, g |^2 ) \, \phi) \, d{v} dx dt = \\ \nonumber
&&\qquad
= \Bigl[ \, \int\!\!\!\int_{\OO}  \beta(g(t,.)) \, \phi \, dx d{v} \, \Bigr]_{t_0}^{t_1}
+ \int_{t_0}^{t_1}\!\!\int\!\!\!\int_{\Sigma} \beta(\gamma \, g) \, \phi \,\,  n(x) \cdot {v} \,\,
d{v} d\sigma_{\! x} dt,  
\eear
holds for all ${t_0}, \,{t_1} \in [0,T]$, all $\beta \in W^{2,\infty}_{loc}(\R)$ and all test functions
$\phi \in \DD([0,T]Ê\times \bar \OO)$. 
\end{theo}
 
\begin{rem} A fundamental point, which is a consequence of the Green formula (\ref{3.10}), is
the possibility of renormalizing the trace function, i.e.
\beqn\label{3.11}
\gamma \, \beta(g) = \beta (\gamma \, g)
\eeqn
for all $\beta \in W^{2,\infty}(\R)$. More generally, (\ref{3.11}) holds as soon as $\gamma \, \beta(g)$
is defined. This is the property that will allow us to define the trace of a renormalized solution.
\end{rem}

\begin{theo}\label{theoTraceLp}
{\bf (The case $p \in [1,\infty)$)}. 
Let $g \in L^\infty(0,T;L^p_{loc}(\bar \OO))$
be a solution of equation (\ref{3.2})-(\ref{(3.4)})-(\ref{(3.6)}).
There exists $\gamma g$ defined on $(0,T) \times \Sigma$ and for every $t \in [0,T]$ there exists $\gamma_t g \in L^p(\OO)$ such that
\beqn\label{3.12}
\gamma_t g \in C([0,T];L^1_{loc}(\OO)) 
\quad \hbox{ and } \quad
\gamma \, g \in L^1_{loc}\bigl([0,T]Ê\times \Sigma, d\lambda_2 \bigr),
\eeqn
and satisfy the Green formula (\ref{3.10}) for every $t_0, \, t_1 \in [0,T]$, 
every $\beta \in \BB_1$ and every test functions  $\phi \in \DD([0,T]Ê\times \bar \OO)$, 
as well as  for every $t_0, \, t_1 \in [0,T]$, every $\beta \in \BB_2$ and every test functions
$\phi \in \DD_0([0,T]Ê\times \bar \OO)$, the space of functions $\phi \in \DD([0,T]Ê\times \bar
\OO)$ such that $\phi = 0$ on $(0,T)Ê\times \Sigma_0$. 
\end{theo}

\begin{theo}\label{theoTraceRenom}
{\bf (The renormalized case)}.  Let $g \in L((0,T) \times \OO)$ satisfy the bound condition (\ref{3.7}) and
the  equation (\ref{3.8}). Then there exists  $\gamma g \in L([0,T]Ê\times \Sigma)$ and for every $t \in [0,T]$ there exists $\gamma_t g \in
C([0,T];L(\OO))$ which satisfy the Green formula (\ref{3.10}) for all ${t_0}, \,{t_1} \in [0,T]$, 
all $\beta \in \BB_1$ and all test functions $\phi \in \DD([0,T]Ê\times \bar \OO)$. 
Furthermore, if (\ref{3.8}) makes sense for at least one function $\beta$ such that $\beta(s) \nearrow
\infty$ when $s \nearrow \infty$, then $\gamma_t g \in L^0(\OO)$ for any $t \in [0,T]$ and
$\gamma g \in L^0([0,T]Ê\times \Sigma)$. 
\end{theo}

%%%%%%%%%%%%%%%%%%%%%%%%%%%%%%%%%%%%%%%%%%%%%%%

\subsection{Proof of the trace theorems}

We begin with some notations. 
For a given real $R > 0$, we define $B_R = \{ y \in \R^N \ / \ |y| < R \}$, $\Omega_R = \Omega
\cap B_R$, $\OO_R = \Omega_R \times B_R$ and $\Sigma_R = (\partial\Omega \cap
B_R) \times B_R$. We also denote by $L^{a,b}_R$
the space $L^a(0,T;L^b(\OO_R))$ or $L^a(0,T;L^b(\Omega_R))$, and $L^{a,b}_{loc}$ the space
$L^a(0,T;L^b_{loc}(\bar \OO))$ or $L^a(0,T;L^b_{loc}(\bar\Omega))$.

\smallskip
\noindent{\sl Proof of Theorem~\ref{theoTraceLinfty}. First step: a priori bounds.} 
In this step we assume that $g$ is a solution of (\ref{3.1}) and is ``smooth". Precisely,
$g \in W^{1,1}\bigl( 0,T; W^{1,\infty}(\Omega; W^{2,\infty}(\R^N) )
\bigr)$, in such a way that the Green formula (\ref{3.10}) holds. The trace $\gamma g$ in (\ref{3.10}) is defined
thanks to the usual trace theorem in the Sobolev spaces.
We shall prove two a priori bounds on $g$. 
Let us define $\beta \in W^{2,\infty}_{loc}(\R)$ by 
$\beta (s) = \cases{ |s| - 1/2 & if $|s| \ge 1$ \cr s^2/2 & if $|s| \le 1$ \cr}$ so that
$\beta' (s) = \cases{ 1 & if $s \ge 1$ \cr s & if $|s| \le 1$ \cr - 1 & if $s \le -1$ \cr}$
and $\beta'' (s) = \cases{ 0 & if $|s| \ge 1$ \cr 1 & if $|s| \le 1$ \cr}$, and thus $\beta \in
\BB_1$. Fix $R>0$ and consider $\chi \in \DD(\bar \OO)$ such that
$0 \le \chi \le 1$, $\chi = 1$ on $\OO_R$ and supp$\, \chi \subset \bar \OO_{R+1}$. We set
$\phi = \chi \,\, n(x) \cdot {v}$. The Green formula (\ref{3.10}) gives
\bean
&&\int_0^T\!\!\!\int\!\!\!\int_{\Sigma} \beta(\gamma \, g) \, \chi \,\,  (n(x) \cdot {v})^2 \,\,
d{v} d\sigma_{\! x} dt =  
- \Bigl[ \, \int\!\!\!\int_{\OO}  \beta(g(t,.)) \, \phi \, dx d{v} \, \Bigr]_0^T \\
&&\qquad\qquad
+\int_0^T \!\!\!\int\!\!\!\int_{\OO} \bigl( \beta(g) \, \Lambda^\star_E \phi 
+ (\beta'(g) \, G - \nu \, \beta''(g) \, |\nabla_{v} \, g |^2 ) \, \phi) \, d{v} dx dt .  
\eean
We deduce from it a first a priori bound: there are some constants $\gamma_R$ and $C_R$ such
that
\bear\nonumber
&&\gamma_R \int_0^T\!\!\!\int\!\!\!\int_{\Sigma_R} |\gamma \, g| \,\,  (n(x) \cdot {v})^2
\,\,d{v} d\sigma_{\! x} dt \le 
\int_0^T\!\!\!\int\!\!\!\int_{\Sigma_R} \beta(\gamma \, g) \,\,  (n(x) \cdot {v})^2
\,\,d{v} d\sigma_{\! x} dt \\ \label{5.1}
&&\le C_R \int_0^T \!\!\!\int\!\!\!\int_{\OO_{R+1}} 
\bigl( g^2 \, (1 + |E|) + |G| + \nu \, |\nabla_{v} \, g |^2 \bigr) \, d{v} dx dt \\ \nonumber
&&\qquad\qquad
+ C_R \int\!\!\!\int_{\OO_{R+1}}  \bigl( g^2(0,.) + g^2(T,.) \bigr) \, dx d{v},  
\eear
where we have used the fact that for $u \in L^\infty(Y_R)$ with $Y_R = \OO_R$ or $\Sigma_R$ there holds
$$
\gamma_R \int_{Y_R}Ê| u | \ \le \ \int_{Y_R} \beta(u) \ \le \ \gamma_R^{-1} \int_{Y_R} u^2.
$$

Let $K \subset \OO$ be a compact set and consider $\phi \in \DD(\OO)$ such that
$0 \le \phi \le 1$, $\phi = 1$ on $K$ and $R > 0$ such that supp$\, \phi \subset \OO_R$. 
We fix $t_0 \in [0,T]$. The Green formula (\ref{3.10}) implies
\bear\label{5.2}
&&\int\!\!\!\int_{\OO}  \beta(g(t_1,.)) \, \phi \, dx d{v}
= \int\!\!\!\int_{\OO}  \beta(g(t_0,.)) \, \phi \, dx d{v} \\ \nonumber
&&\qquad\qquad +\int_{t_0}^{t_1}\!\int\!\!\!\int_{\OO} \bigl( \beta(g) \, \Lambda^\star_E \phi 
+ (\beta'(g) \, G - \nu \, \beta''(g) \, |\nabla_{v} \, g |^2 ) \, \phi) \, d{v} dx dt,  
\eear
and we get a second a priori bound
\bear\label{5.3}
\gamma_R \int\!\!\!\int_K  |g| (t_1,.) \, dx d{v} &\le& C_R \int\!\!\!\int_{\OO_R}  g^2(t_0,.) \, dx d{v} \\ \nonumber
&+& C_R \int_0^T \!\!\!\int\!\!\!\int_{\OO_R} 
\bigl( g^2 \, (1 + |E|) + |G| + \nu \, |\nabla_{v} \, g |^2 \bigr) \, d{v} dx dt. 
\eear

\medskip
\noindent{\sl Second step: regularization and passing to the limit.} 
Let us now consider a function $g$ which satisfies the assumptions of Theorem~\ref{theoTraceLinfty}.
We define the mollifier $\rho_k$  by
$$
\rho_k(z) = k^N \, \rho(k \, z) \ge 0, \quad k \in \N^\star, \quad \rho \in \DD(\R^N), 
\quad \hbox{supp} \, \rho \subset B_1, \quad \int_{\R^N} \rho(z) \, dz = 1,
$$
and we introduce the regularized functions $g_k = g \star_{x,k} \rho_k *_{v} \rho_k$,
where $*$ stands for the usual convolution and $\star_{x,k}$ for the convolution-translation defined
by
$$
(u \star_{x,k} h_k)(x) = [\tau_{2 \, n(x) /k} (u * h_k)] (x) 
= \int_{\R^N} u(y) \, h_k (x - {2 \over k} \, n(x) - y ) \, dy,
$$
for all $u \in L^1_{loc}(\bar\Omega)$ and $h_k \in L^1(\R^N)$ with supp$\, h_k \subset B_{1/k}$.

\begin{lem}\label{Lemma51}
With this notation one has
$g_k \in W^{1,1}\bigl( 0,T; W^{1,\infty}(\Omega; W^{2,\infty}(\R^N) ) \bigr)$  and
$$
\Lambda_E g_k = G_k \ \ \hbox{ in } \ \  \DD'((0,T) \times \OO), 
$$
with $G_k \in L^1_{loc}((0,T) \times \bar\OO)$ for all $k \in \N$. Moreover, the sequences
$(g_k)$ and $(G_k)$ satisfy
\beqn  \label{5.5}
 \left\{
 \begin{array}{l}
 \displaystyle{
(g_k) \hbox{ is bounded in } L^\infty((0,T) \times \OO),
\quad g_k \longrightarrow g \ \hbox{ a.e. in } (0,T) \times \OO, }\vspace{0.3cm} \\
 \displaystyle{ \nabla_{v} g_k \longrightarrow \nabla_{v} g \ \hbox{ in } L^2_{loc}([0,T] \times \bar\OO)
\quad \hbox{ and } \quad
G_k \longrightarrow G \ \hbox{ in } L^1_{loc}([0,T] \times \bar \OO)}.
 \end{array}
 \right.
\eeqn
\end{lem}
 
\smallskip\noindent
The proof of  Lemma~\ref{Lemma51} is similar to the proof of \cite[Lemma 1]{53} and of \cite[Lemma
II.1]{36} to which we refer.

\smallskip
From Lemma~\ref{Lemma51} we have that for all $k,\ell \in \N^\star$ the difference $g_k - g_\ell$
belongs to
$W^{1,1}\bigl(0,T;$ \break $W^{1,\infty}(\Omega; W^{2,\infty}(\R^N) ) \bigr)$ and is a solution of
$$ 
\Lambda_E (g_k - g_\ell) = G_k - G_\ell \ \ \hbox{ in } \ \  \DD'((0,T) \times \OO).
$$
We know, thanks to (\ref{5.5}), that $g_k(t,.)$ converges to $g(t,.)$ in $L^2_{loc}(\bar\OO)$ for a.e. 
$t \in [0,T]$; we fix $t_0$ such that $g_k(t_0,.) \to g(t_0,.)$. Moreover, up to
a choice for the continuous representation of $g_k$, we can assume that $g_k \in
C([0,T],L^1_{loc}(\bar\OO))$. Therefore, the estimate (\ref{5.2}) applied to $g_k - g_\ell$ in $t_0$ and
the convergence (\ref{5.5}) imply that for all compact sets $K \subset \OO$ we have
\beqn\label{5.7}
\sup_{t \in [0,T]} \| (g_k - g_\ell)(t,.) \|_{L^1(K)}  \mathop{\longrightarrow}_{k,\ell \to
+\infty} 0.
\eeqn
We deduce from this, that there exists, for any time $t \in [0,T]$, a function $\gamma_t g$ such
that $g_k(t,.)$ converges to $\gamma_t g$ in $C([0,T]; L^1_{loc}(\OO))$; in particular,
$$
g(t,x,{v}) = \gamma_t g (x,{v}) \hbox{ for a.e. } (t,x,{v}) \in (0,T) \times \OO.
$$
Thus, we also have $g_k(t,.) = (\gamma_t \, g) \star_{x,k} \rho_k *_{v} \rho_k$ a.e. in $(0,T) \times \OO$, and
since these two functions are continuous, the equality holds for all $(t,x,{v})
\in [0,T] \times \bar\OO$ and $k \in \N^\star$, so that $g_k(t,.) \to \gamma_t g$ in
$L^2_{loc}(\bar \OO)$ for all $t \in [0,T]$. 

\smallskip
Using now the estimate (\ref{5.1}), applied to $g_k - g_\ell$, and the convergence (\ref{5.5}) and (\ref{5.7}) we
get that
$$
\int_0^T\!\!\!\int\!\!\!\int_{\Sigma_R} |\gamma g_k - \gamma g_\ell| \, 
(n(x) \cdot {v})^2 \, d{v} d\sigma_{\! x} dt
\mathop{\longrightarrow}_{k,\ell \to +\infty} 0, 
$$
for all $R > 0$. We deduce that there exists a function $\gamma g \in L^1_{loc}([0,T] \times
\Sigma,(n(x) \cdot {v})^2 \, d{v} d\sigma_{\! x} dt)$, which is the limit of $\gamma
g_k$ in this space. Moreover, since
$ \| \gamma g_k \|_{L^\infty} \le \| g_k \|_{L^\infty}$ is bounded, we have
$\gamma g \in L^\infty((0,T) \times \OO)$.

Finally, we obtain the Green formula (\ref{3.10})  writing it first for $g_k$ and then passing to the
limit $k \to \infty$ thanks to the convergence previously obtained. 
Uniqueness of the trace function follows from the Green formula. \qed

\medskip
\noindent{\sl Proof of Theorem~\ref{theoTraceRenom}. } The proof is based on Theorem~\ref{theoTraceLinfty} and on a monotony argument.
This is exactly the same as the one presented in \cite{54} in the case of Vlasov equation.
Let $(\beta_M)_{M \ge 1}$ be a sequence of
odd functions of $\BB_1$ such that
$$
\beta_M(s) = \cases{
s & if $s \in [0,M]$ \cr
M+1/2 & if $ s \ge M+1$, \cr}
$$
and $|\beta_M(s)| \le |s|$ for all $s \in \R$. 
The function $\alpha_M(s) := \beta_M(\beta^{-1}_{M+1}(s))$, with the convention $\alpha_M(s) = M + 1/2$ if $s \ge M + 3/2$,
is well defined, odd and also belongs to $\BB_1$. 
We will construct the trace function $\gamma g$ as the limit of $(\gamma \beta_M(g))$ when $M
\to \infty$, that one being defined thanks to Theorem~\ref{theoTraceLinfty}.
Indeed, the condition (\ref{(3.6)}) implies that
$$
\nabla_{v} g \, {\bf 1}_{|g| \le M+1} \in L^2_{loc}([0,T] \times \bar \OO),
$$
and then $\nabla_{v} \beta_M(g) = \beta'_M(g) \, \nabla_{v} g \in  L^2_{loc}([0,T] \times \bar \OO)$
in such a way that $\beta_M (g)$ satisfies the assumption on Theorem~\ref{theoTraceLinfty}.
We define
$\Gamma_M^{(\pm)} = \{ (t,x,{v}) \in (0,T) \times \Sigma, \pm \gamma \beta_M(g) (t,x,{v}) > 0 \}$
and
$\Gamma_M^{(0)} = \{ (t,x,{v}) \in (0,T) \times \Sigma, \gamma \beta_M(g) (t,x,{v}) = 0 \}$. 
Thanks to the definition of $\alpha_M$ and the renormalization property (\ref{3.11}) of the
trace, one has $\gamma \, \beta_M(g) = \gamma \, \alpha_M(\beta_{M+1}(g)) = \alpha_M(\gamma
\,\beta_{M+1}(g))$. We deduce that, up to a set of measure zero,
$$ 
\Gamma_M^{(+)} = \Gamma_1^{(+)}, \quad \Gamma_M^{(-)} = \Gamma_1^{(-)} \quad
\hbox{and} \quad \Gamma_M^{(0)} = \Gamma_1^{(0)} \quad \hbox{for all } M \ge 1.
$$
Therefore the sequence $(\gamma \, \beta_M(g))_{M \ge 1}$ is increasing on $\Gamma_1^{(+)}$ and
decreasing on $\Gamma_1^{(-)}$. This implies that $\gamma \, \beta_M(g)$
converges a.e. to a limit denoted by $\gamma g$ which belongs to
$L([0,T] \times \Sigma)$. Obviously, if (\ref{3.8}) holds for one function $\beta$ such that $\beta(s)
\nearrow +\infty$ when $s \nearrow \pm\infty$, then $\beta(\gamma g) \in L^1((0,T) \times
\Sigma,d\lambda_2)$ and $\gamma g \in L^0((0,T) \times \Sigma)$.
In order to establish the Green formula (\ref{3.10})  we fix $\beta \in \BB_1 $ and 
$\phi \in \DD((0,T] \times \bar \OO)$. We write the Green formula for the function
$\beta(\beta_M(g))$, and using the fact that $\gamma \bigl[ \beta \circ \beta_M(g) \bigr] =
\beta(\gamma \beta_M(g))$, we find
\bean
&&\int_{0}^{T}\!\int\!\!\!\int_{\OO} \bigl( \beta \circ \beta_M (g) \, 
({\partial \phi \over \partial t}  +  {v} \cdot \nabla_x  \phi + E \cdot \nabla_{v} \phi ) 
+ (\beta \circ \beta_M)'(g) \, G \, \phi\bigr) \, d{v} dx dt = \\
&&\qquad\qquad=\int_{0}^{T}\!\int\!\!\!\int_{\Sigma} \beta(\gamma \,\beta_M(g)) \, \phi \, n(x) \cdot {v} \, d{v}
d\sigma_{\! x} ds.
\eean
We get (\ref{3.10})  by letting $M \to \infty$ and noticing that
$\beta \circ \beta_M(s) \to \beta(s)$ for all $s \in \R$. \qed

\begin{rem} Theorem~\ref{theoTraceLp} is now a quite simple consequence  of Theorem~\ref{theoTraceRenom} using the
a priori bounds stated in the proof of Theorem~\ref{theoTraceLinfty}. Let  us emphasize that with the additional
assumption (\ref{(3.5)}) in hands, it is  possible to give a direct proof of Theorem~\ref{theoTraceLp} (following the
proof of Theorem~\ref{theoTraceLinfty}) instead of passing through the renormalization step. See \cite{53} for details.
\end{rem}

\smallskip
\noindent{\sl Proof of Theorem~\ref{theoTraceLp}.} For all $\beta \in \BB_1$ it is clear that $\beta(g) \in
L^\infty$, $\nabla_{v} \beta(g) \in L^2$ and that $\beta(g)$ is solution of (\ref{3.8}) using Lemma~\ref{Lemma33} below
(we just have to multiply equation (\ref{5.9}), in the case $\mu \equiv 0$, by $\beta'(g_k)$ and to pass to
the limit $k \to \infty$). Thanks to Theorem~\ref{theoTraceRenom}, we already know that $g$ has a trace
$\gamma_t g \in L(\OO)$ and $\gamma g \in L((0,T) \times \OO)$ which satisfies the Green formula
(\ref{3.10})  for all $\beta \in \BB_1$ and $\phi \in \DD([0,T] \times \bar \OO)$. 
We just have to prove that $\gamma g$ and $\gamma_t g$ belong to the appropriate space.
On one hand, for all $\beta \in \BB_1$ such that $|\beta(s)| \le |s|$ one has
$$
\| \beta(\gamma_t g) \|_{L^p_R} \le \sup_k \sup_{[0,T]} \| \beta(g_k(t,.)) \|_{L^p_R}
\le \sup_{[0,T]} \| g_k(t,.) \|_{L^p_R} \le  \| g \|_{L^{\infty,p}_R},
$$
and thus, choosing $\beta = \beta_M$, defined in the proof of Theorem~\ref{theoTraceRenom}, one gets, passing to the
limit $M \to \infty$,
$$
\sup_{[0,T]} \|\gamma_t g \|_{L^p_R} \le   \| g \|_{L^{\infty,p}_R} < \infty.
$$
In the same way and using (\ref{5.1}), we show that
$$
\| \gamma g \|_{L^1([0,T] \times \Sigma_R, d\lambda_2)} < \infty.
$$
We still have to prove that $\gamma_t g \in C([0,T],L^1_{loc}(\bar\OO))$, which is an immediate
consequence of the following Lemma. \qed

\begin{lem}\label{Lemma5.2}
Let $(u_n)$ be a bounded sequence of $L^1_{loc}(\OO)$ such that
$\beta(u_n) \rightharpoonup \beta(u)$ in $\bigl( C_c(\OO) \bigr)'$ for all $\beta \in \BB_2$.
Then $u_n \to u$ in $L^1_{loc}(\OO)$.
\end{lem}

\smallskip\noindent
{\sl Proof of Lemma~\ref{Lemma5.2}.} We fix $j: \R \to \R$ a non-negative function of class $C^2$,
strictly convex on the interval $[-M,M]$ and such that $j''(t) = 0$ for all $t \notin
[-M,M]$; in particular
$j \in \BB_2$. We also consider $\chi \in C_c(\OO)$ such that $0 \le \chi \le 1$.
By assumption
\beqn\label{5.10}
\int_\OO j(u_n) \, \chi \to \int_\OO j(u) \, \chi 
\eeqn
and by convexity of $j$ one also has
\beqn\label{5.11}
\liminf_{n \to \infty} \int_\OO j \bigl( {u_n + u \over 2} \bigr) \, \chi
\ge \int_\OO j(u) \, \chi
\quad \hbox{ since } \quad
{u_n + u \over 2}Ê\rightharpoonup u \ \hbox{ in } \ \bigl( C_c(\OO) \bigr)'.
\eeqn
Remarking that
\beqn\label{5.12}
{1 \over 2} \, j(t) + {1 \over 2} \, j(s) - j \bigl( {t+s \over 2} \bigr) \ge 0 
\qquad \forall t , \, s \in \R,
\eeqn
we deduce from (\ref{5.10}) and (\ref{5.11}) that 
\beqn\label{5.13}
\int_\OO \bigl[ {1 \over 2} \, j(u_n) + {1 \over 2} \, j(u)
- j \bigl( {u_n+u \over 2} \bigr) \bigr] \, \chi \to 0.
\eeqn
From the fact that in (\ref{5.12}) the inequality is strict whenever $t , \, s \in
[-M,M]$ and $t \not=s$, we obtain from (\ref{5.13}) that there exists a subsequence $(u_{n_k})$ such
that $u_{n_k} \to u$ a.e. on  supp$\, \chi \cap [ \, |u| < M]$.
The preceding argument being valuable for arbitrary $M$ and $\chi$, we obtain, by a diagonal
process, a subsequence of $(u_n)$, still denoted by $(u_{n_k})$, such that
$u_{n_k} \to u$ a.e. in $\OO$.

We now set $j_\pm(s) = s_\pm$. We first remark that we can write $j_\pm = j_{\pm,1} + j_{\pm,2}$
with $j_{\pm,1} \in \BB_2$ and $j_{\pm,2} \in W^{2,\infty}(\R)$ in such a way that
$$
\int_\OO j_\pm (u_{n_k}) \, \chi \to \int_\OO j_\pm (u) \, \chi.
$$
On the other hand, the elementary inequality $\bigl| \,  b - |a - b| \, \bigr| \le a$ $ \
\forall a, \, b \ge 0$ and the dominated convergence Theorem imply
$j_\pm(u_{n_k}) - |j_\pm(u_{n_k}) - j_\pm(u)| \to j_\pm(u)$ in $L^1_{loc}(\OO)$. It follows that
$$
\limsup_{k \to \infty} \int_\OO \bigl| \, j_\pm(u_{n_k}) - j_\pm(u) \, \bigr| \, \chi 
= \int_\OO j_\pm(u) \, \chi - \lim_{k \to \infty} \int_\OO j_\pm(u_{n_k}) \, \chi = 0.
$$
We conclude that $u_{n_k} = j_+( u_{n_k}) - j_-( u_{n_k}) \to j_+(u) - j_-(u) = u$ strongly in
$L^1_{loc} (\OO)$ and that, in fact, it is the whole sequence $(u_n)$ which converges.
\qed

%%%%%%%%%%%%%%%%%%%%%%%%%%%%%%%%%%%%%%%%%%%%%%%%
\section{Renormalized convergence for the trace functions sequence.}
\label{sec:GalStab}
\setcounter{equation}{0}
\setcounter{theo}{0}
%%%%%%%%%%%%%%%%%%%%%%%%%%%%%%%%%%%%%%%%%%%%%%%%

%!!!!!!!!!!!!!!!!!!!!!!!!!!!!!!!!!!!!!!!!!!!!!!!!!!!!!!!!!!!!!!!!!!!!!!!!!!!!!!!!!!!!!!!
%         Renormalized convergence for trace functions sequence
%!!!!!!!!!!!!!!!!!!!!!!!!!!!!!!!!!!!!!!!!!!!!!!!!!!!!!!!!!!!!!!!!!!!!!!!!!!!!!!!!!!!!!!!

\medskip
We present now a quite general stability result in both the interior and up to the boundary for a
sequence of renormalized solutions to the Vlasov-Fokker-Planck equation in a domain. This
will be a key argument in the proof of Theorem~\ref{ThStab}.
In some sense, this result says that renormalized convergence, as well as the a.e. convergence, can be
propagated from the interior to the boundary. Notice that it is not clear that a similar result holds for the $L^1$-weak convergence. 

\begin{theo}\label{th:rcvgeV}
Define $\BB_3$ as the class of functions of $W^{1,\infty}_{loc}(\R)$ such that $|\beta'(s)| \, (1+s)^{-1} \in L^\infty(\R)$. 
Consider three sequences $(g_n)$, $(E_n)$ and $(G_n)$, with $G_n = G^+_n - G^-_n$, $G^\pm_n \ge 0$, which satisfy 
for any renormalizing sequence $(\alpha_M)$ in $\BB_3$ and for any $\beta \in \BB_3$ the convergence assumptions
\bear
\label{stabA1}
&& g_n \wto g  \,\, \hbox{weakly in} \,\, L^\infty(0,T: L^1 (\OO)), \\
\label{stabA2}
&& E_n \to E \,\, \hbox{strongly in} \,\,  L^1((0,T) \times \Omega), \hbox{uniformly bounded in } L^1(0,T; W^{1,1}(\Omega)), 
\\ \label{stabA3}
&& \alpha_M'(g_n) \, G^\pm_n \rightharpoonup \bar G^\pm_M \,\,  \hbox{weakly in } L^1((0,T) \times \OO_R), \\Ê\nonumber 
&&\quad\hbox{with} \,\, 
\bar G^\pm_M \nearrow G^\pm \,\,  \hbox{ a.e. and }Ê\, \beta'(g) \, G^\pm \in L^1((0,T) \times \OO_R), 
\eear
as well as the renormalized Vlasov equation 
\beqn\label{stabA4}
\Lambda_{E_n} \, \beta( g_n ) = \beta'(g_n) \, G_n \  \hbox{ in }  \ \DD'((0,T) \times \OO),  
\eeqn
for which each term clearly makes sense thanks to (\ref{stabA1})--(\ref{stabA3}).
Then $g \in L^\infty(0,T; L^1(\OO))$ is a solution of
\beqn\label{stabA6}
\Lambda_E \, \beta( g ) = \beta'(g) \, G 
\  \hbox{ in } \ \DD'((0,T) \times \OO), \quad G = G^+-G^-, 
\eeqn
for any $\beta \in \BB_3$. Furthermore, the traces $\gamma g_n$ and $\gamma g$ defined thanks to
the Theorem~\ref{theoTraceRenom} satisfy
\beqn\label{stabA5}
\gamma g_n \rto \gamma g \quad\hbox{in the renormalized sense.}
\eeqn
\end{theo}

\noindent{\sl Proof of Theorem~\ref{th:rcvgeV}. } The proof is essentially the same as Step 2 in the proof of \cite[Proposition 5]{54} and 
as the proof of Theorem~\ref{FTeqStab}. Nevertheless, for the sake of completeness, we sketch the main arguments.

\smallskip\noindent
{\sl Step 1. }ÊUp to the extraction of a  subsequence, we have $\displaystyle{g_n \rto g}$ thanks to (\ref{stabA1}) and Lemma~\ref{BitToRenorm}, and there exists $\eta \in L((0,T) \times \Sigma)$ such that $\displaystyle{\gamma g_n \rto \eta}$ thanks to  Proposition~\ref{propWW1}. More precisely, there exists two sequences $(\bar \alpha_M)$ and $(\bar \gamma_M)$ and  such that 
\bear\label{stabA10}
&&\alpha_M(g_{n}) \wto \bar \alpha_M \quad \sigma(L^\infty,L^1) \, \star
\quad\hbox{ and }\quad \bar \alpha_M \nearrow g \quad \hbox{ a.e.}, \\
\label{stabA11}
&&\alpha_M(\gamma g_{n}) \wto \bar \gamma_M  \quad \sigma(L^\infty,L^1) \, \star
\quad\hbox{ and }\quad \bar \gamma_M  \nearrow \eta \quad \hbox{ a.e.}.
\eear
The  Green formula (\ref{3.10}) associated to the equation (\ref{stabA4}) with $\beta = \alpha_M$ implies
$$
\int_{0}^T\!\!\int\!\!\!\int_{\OO} \bigl( \alpha_M(g_{n})   \, \Lambda^*_E \varphi +  \alpha_M'(g_n) \, G_n  \varphi \bigr) \, d{v} dx dt 
=  \int_0^T\!\!\int\!\!\!\int_{\Sigma} \alpha_M(\gamma \, g_n) \, \varphi \,\,  n(x) \cdot {v} \,\,
d{v} d\sigma_{\! x} dt, 
$$
for any $\varphi \in \DD((0,T) \times \bar\OO)$.  Passing to the limit  $M \to \infty$ with the help of (\ref{stabA10}), (\ref{stabA2}), (\ref{stabA3}) in the above identity, we obtain  
\beqn\label{stabA12}
\Lambda_E \, \bar\alpha_M = \bar G_M := \bar G_M^+ - \bar G_M^- \quad\hbox{in}\quad \DD'((0,T) \times \OO),
\eeqn
and $\gamma \bar \alpha_M = \bar\gamma_M$ thanks to the trace Theorem~\ref{theoTraceLinfty} and the convergence (\ref{stabA11}). 

\smallskip\noindent
{\sl Step 2. } For a given function $\beta \in \BB_3 \cap L^\infty$, we write the renormalized Green formula (\ref{3.10}) associated to the equation (\ref{stabA12}) as
\beqn\label{stabA13}
\int_{0}^T\!\!\int\!\!\!\int_{\OO} \bigl( \beta(\bar \alpha_M)   \, \Lambda^*_E \varphi +  \beta'(\bar \alpha_M) \, \bar G_M  \varphi \bigr) \, d{v} dx dt  =  \int_0^T\!\!\int\!\!\!\int_{\Sigma} \beta(\gamma_M) \, \varphi \,\,  n(x) \cdot {v} \,\,
d{v} d\sigma_{\! x} dt, 
\eeqn
for any $\varphi \in \DD((0,T) \times \bar\OO)$. 
Using that $(\bar \alpha_M)$, $(G^\pm_M)$ and $(\bar \gamma_M)$ are a.e. increasing sequences we have 
\beqn\label{stabA14}
\beta(\bar\alpha_M) \nearrow \beta(g), \quad \beta'(\bar\alpha_M) \, \bar G_M^\pm \nearrow \beta'(g) \, G^\pm \quad\hbox{in}\quad L^1((0,T) \times \OO)
\eeqn
as well as 
\beqn\label{stabA15}
\beta(\bar\gamma_M) \nearrow \beta(\eta) \,\, \hbox{a.e. and uniformaly bounded in } \, L^\infty((0,T) \times \OO). 
\eeqn
Passing to the limit in (\ref{stabA13}) with the help of (\ref{stabA14}) and (\ref{stabA15}) we obtain
$$
\int_{0}^T\!\!\int\!\!\!\int_{\OO} \bigl( \beta(g)  \, \Lambda^*_E \varphi +  \beta'(g) \, G \,  \varphi \bigr) \, d{v} dx dt  =  \int_0^T\!\!\int\!\!\!\int_{\Sigma} \beta(\eta) \, \varphi \,\,  n(x) \cdot {v} \,\,
d{v} d\sigma_{\! x} dt, 
$$
which precisely means that $\eta = \gamma g$. We conclude by gathering that information with (\ref{stabA11}). \qed

\begin{theo}\label{th:rcvgeVFP}
Consider three sequences $(g_n)$, $(E_n)$ and $(G_n)$ which satisfy, for
all $\beta \in \BB_4$ the class of functions of $W^{2,\infty}_{loc}(\R)$ such that $|\beta'(s)| \,
(1+s)^{-1} \in L^\infty(\R)$ and $|\beta''(s)| \, (1+s)^{-2} \in L^\infty(\R)$,
\bear
\label{3.13}
&& g_n \rightarrow g  \hbox{ strongly in } L^1((0,T) \times \OO) \hbox{ and is uniformly
bounded in } L^\infty(0,T; L^1(\OO)), \\
\label{3.14}
&& E_n \rightharpoonup E \quad \hbox{ weakly in } L^1(0,T; W^{1,1}_{loc}(\bar\OO)),
\\ \label{3.15}
&& \beta'(g_n) \, G_n \rightharpoonup \beta'(g) \, G \quad \hbox{weakly in } 
L^1((0,T) \times \OO_R), \   \forall R \ge 0, 
\\ \label{3.16}
&&  \!\!\int_0^T \int_\OO  {|\nabla_{v} g_n|^2 \over 1 + g_n} \, d{v} dx dt \le C_T,
\eear
as well as the renormalized Vlasov-Fokker-Planck equation 
\beqn\label{3.17}
\Lambda_{E_n} \, \beta( g_n ) = \beta'(g_n) \, G_n - \nu \, \beta''(g_n) \, | \nabla_{v} g_n|^2
\quad \hbox{ in } \DD'((0,T) \times \OO),
\eeqn
for which each term makes sense thanks to (\ref{3.13})--(\ref{3.16}).
Then $g \in L^\infty(0,T; L^1(\OO))$ is a solution of
\beqn\label{3.18}
\Lambda_E \, \beta( g ) = \beta'(g) \, G - \nu \, \beta''(g) \, | \nabla_{v} g|^2
\quad \hbox{ in } \DD'((0,T) \times \OO)
\eeqn
for all $\beta \in \BB_4$. Furthermore, the traces $\gamma g_n$ and $\gamma g$ defined thanks to
the Theorem~\ref{theoTraceRenom} satisfy
\beqn\label{3.19}
\gamma g_n \rto \gamma g \quad\hbox{in the renormalized sense, and}\quad
\gamma_+ g_n \rightarrow \gamma_+ g \quad a.e.
\eeqn
\end{theo}

\noindent
We shall need the following auxiliary results in the proof of Theorem~\ref{th:rcvgeVFP}.

\begin{lem}\label{Lemma31}
Let $(u_n)$ be a bounded sequence of $L^2(Y)$ such that $u_n \rightharpoonup u$ weakly in
$L^2(Y)$. Then, there exists $\mu \in (C_c(Y))'$, a non-negative measure, such that, up to
the extraction of a subsequence,
$$
|u_n|^2 \rightharpoonup |u|^2 + \mu 
\qquad \hbox{ weakly in } \quad
(C_c(Y))'.
$$
\end{lem}

\begin{lem}\label{Lemma32} For any $\theta \in (0,1)$ and $M \in (0,\infty)$ we set
$$
\Phi (s) = \Phi_{M,\theta} (s) := \cases{ {1 /\theta} \, (e^{\theta \, s} - 1) & if $s \le M$ \cr  
(s - M) \, e^{\theta \, M} + {1 /\theta} \, (e^{\theta \, M} -1)& if $s \ge M,$}
$$ 
and $\beta(s) := \beta_1(s) = \log (1+s)$. Then
$$ \left\{
 \begin{array}{l}
 \displaystyle{
 \Phi'(s) \ge 1, \qquad
\Phi \circ \beta (s) \nearrow s \hbox{ when } M \nearrow \infty, \ \theta \nearrow 1, }\vspace{0.3cm} \\
 \displaystyle{  \hbox{and} \qquad 0 \le - (\Phi \circ \beta)''(s) \le {1-\theta + e^{(\theta-1) \, M} \over 1 + s}}
 \end{array}
 \right. \qquad\quad \forall s \ge 0.
$$
\end{lem}

\begin{lem}\label{Lemma33} Let $g \in L^\infty(0,T;L^p_{loc}(\OO))$ be a solution to the
Vlasov-Fokker-Planck equation
\beqn\label{5.9}
\Lambda_E \, g =  G + \mu \ \hbox{ in } \DD'((0,T) \times \OO),
\eeqn
with $E \in L^1(0,T;W^{1,p'}_{loc}(\OO))$, $G \in L^1_{loc}((0,T) \times \OO))$ and
$\mu \in \DD'((0,T) \times \OO)$, $\mu \ge 0$.
For a given mollifer $\rho_k$ in $\R^N$, we set
$$
g_k := g *_t \rho_k *_x \rho_k *_{v} \rho_k
\quad\hbox{ and }\quad
\mu_k := \mu *_t \rho_k *_x \rho_k *_{v} \rho_k.
$$
Then $g_k$ satisfies the Vlasov-Fokker-Planck equation
$$
\Lambda_E \, g_k =  G_k + \mu_k \ \hbox{ in all compact set of } (0,T) \times \OO,
$$
with $G_k \to G$ strongly in $L^1_{loc}([0,T] \times \OO))$.
\end{lem}

\smallskip\noindent
The proof of Lemma~\ref{Lemma31} is classical, the one of Lemma~\ref{Lemma32} is elementary, and we refer to 
\cite{35} for the proof of Lemma~\ref{Lemma33}.

\smallskip\noindent
{\sl Proof of the Theorem~\ref{th:rcvgeVFP}.} {\sl Step 1: Proof of (\ref{3.18}).} This step is inspired from \cite{35} and
it is clear from the theory of renormalized solution \cite{36} that it is enough to prove (\ref{3.18})
only for $\beta(s) := \log (1 +s)$. With the notation $h_n := \beta(g_n)$ and $h = \beta(g)$ we have
$\nabla_{v} h_n = \sqrt{- \beta''(g_n)} \, \nabla_{v} g_n \rightharpoonup 
 \sqrt{- \beta''(g)} \, \nabla_{v} g = \nabla_{v} h$ weakly in $L^2((0,T) \times \OO)$ so that,
thanks to Lemma~\ref{Lemma31}, there is a bounded measure $\mu \ge 0$ such that, up to the extraction of a
subsequence, $ | \nabla_{v} h_n|^2 \rightharpoonup | \nabla_{v} h |^2 + \mu$ weakly in $\DD'([0,T]
\times \bar\OO)$. Passing to the limit $n \to \infty$ in (\ref{3.17}) we get
$$
\Lambda_E \, \beta(g)  = \beta'(g) \, G - \beta''(g) \, | \nabla_{v} g|^2 + \mu
\quad \hbox{ in } \DD'((0,T) \times \OO).
$$
We just point out that 
$$
E_n \, \beta(g_n) \rightharpoonup E \, \beta(g)
\quad\hbox{ weakly in }L^1((0,T) \times \OO),
$$
since $\beta(g_n) \rightarrow \beta(g)$ strongly in $L^2(0,T;L^p(\OO))$ for all $p < \infty$
and $E_n  \rightharpoonup E$ weakly in $L^2(0,T;L^q(\OO))$ for every $q \in [1,N/(N-1))$. 
We prove now that $\mu = 0$ in $(0,T) \times \OO$.

\smallskip\noindent
With the notations introduced in Lemma~\ref{Lemma32} and Lemma~\ref{Lemma33} we have
$$
\Lambda_E \, \Phi ( h_k)  = 
\Phi'(h_k) \, \bigl( \beta'(g) \, G -  \beta''(g) \, | \nabla_{v} g|^2 \bigr) *_{t,x,{v}} \rho_k 
- \Phi''(h_k) \, |\nabla_{v} h_k|^2 + \Phi'(h_k) \, \mu_k.
$$
Using that $\Phi' \ge 1$ (thanks to Lemma~\ref{Lemma32}) and passing to the limit $k \to \infty$ (thanks to  Lemma~\ref{Lemma33}), we get
$$
\Lambda_E \, (\Phi \circ \beta) (g) \ge \Phi'(\beta(g)) \, \beta'(g) \, G 
- (\Phi'(\beta(g)) \, \beta''(g) + \Phi''(\beta(g)) \, (\beta'(g))^2 ) \, | \nabla_{v} g|^2 + \mu
$$
and then
\beqn\label{3.20}
\Lambda_E \, (\Phi \circ \beta) (g) - (\Phi \circ \beta)'(g) \, G 
\ge (\Phi \circ \beta)'' (g) \, | \nabla_{v} g|^2 + \mu
\quad \hbox{ in } \DD'((0,T) \times \OO).
\eeqn
In order to have an estimate of the left hand side we come back to equation (\ref{3.17}), and we
write
$$
\Lambda_{E_n} \, \Phi \circ \beta ( g_n ) = (\Phi \circ \beta)'(g_n) \, G_n 
- (\Phi \circ \beta)''(g_n) \, | \nabla_{v} g_n|^2
\quad \hbox{ in } \DD'((0,T) \times \OO)
$$
since $\Phi \circ \beta \in \BB_4$.
Then, for all $\chi \in \DD((0,T) \times \OO$ such that $0 \le \chi \le 1$ we have
(thanks to Lemma~\ref{Lemma32})
\bean
&&\Bigl| \int_0^T \!\! \int_\OO \bigl( \Phi \circ \beta ( g_n ) \, \Lambda_{E_n} \, \chi
+ (\Phi \circ \beta)'(g_n) \, G_n \, \chi \bigr) \, d{v} dx dt \Bigr| = \\
&&\qquad\qquad
= -\int_0^T \!\! \int_\OO (\Phi \circ \beta)''(g_n) \, | \nabla_{v} g_n|^2 \, \chi \, d{v} dx dt
\\
&&\qquad\qquad
\le [1-\theta + e^{(\theta-1) \, M}] \, 
\int_0^T \!\! \int_\OO {| \nabla_{v} g_n|^2 \over 1 + g_n} \, d{v} dx dt.
\eean
Passing to the limit $n \to \infty$ we get, thanks to (\ref{3.16}),
$$
\Bigl| \int_0^T \!\!\! \int_\OO \bigl( \Phi \circ \beta ( g ) \, \Lambda_E \, \chi
+ (\Phi \circ \beta)'(g) \, G \, \chi \bigr) \, d{v} dx dt \Bigr| 
\le [1-\theta + e^{(\theta-1) \, M}] \, C_T.
$$
Then, coming back to (\ref{3.20}), we have (thanks to Lemma~\ref{Lemma32} again)
\bean
\int_0^T \!\!\! \int_\OO \chi \, d\mu &\le&
- \int_0^T \!\!\! \int_\OO \bigl( \Phi \circ \beta ( g ) \, \Lambda_E \, \chi
+ (\Phi \circ \beta)'(g) \, G \, \chi 
+ (\Phi \circ \beta)''(g) \, | \nabla_{v} g|^2 \bigr) \, d{v} dx dt \\
&\le& 2 \, [1-\theta + e^{(\theta-1) \, M}] \, C_T
\qquad \forall \theta \in [0,1], \ M > 0, 
\eean
and letting $M \to \infty$ and then $\theta \to 1$ we obtain $\mu = 0$ on supp$\, \chi$,
which is precisely saying that $\mu = 0$ on $(0,T) \times \OO$.

\medskip
\noindent{\sl Step 2: Proof of (\ref{3.19}).}
We fix $\phi \in \DD((0,T) \times \bar\OO)$ such that $0 \le \phi \le 1$. By definition of 
$\gamma g_n$ we have
\bean
&& \Bigl| 
\int_0^T\!\int\!\!\!\int_{\Sigma} \Phi \circ \beta(\gamma \, g_n) \, \phi \,\,  n(x) \cdot {v}
\,\, d{v} d\sigma_{\! x} dt \\
&&\qquad\qquad\qquad- \int_0^T \!\! \int_\OO \bigl( \Phi \circ \beta ( g_n ) \, \Lambda_{E_n} \, \chi
+ (\Phi \circ \beta)'(g_n) \, G_n \, \chi \bigr) \, d{v} dx dt  \Bigr| =\\
&&\qquad= \int_0^T \!\! \int_\OO (\Phi \circ \beta)''(g_n) \, | \nabla_{v} g_n|^2 \, \chi \, d{v} dx dt
\le [1-\theta + e^{(\theta-1) \, M}] \, 
\int_0^T \!\! \int_\OO {| \nabla_{v} g_n|^2 \over 1 + g_n} \, d{v} dx dt.\
\eean
We note $\overline{\Phi \circ \beta}$ the $L^1$-weak limit of $\Phi \circ \beta(\gamma \, g_n)$.
Passing to the limit $n \to \infty$ we get
\bean
&&\Bigl| 
\int_0^T\!\int\!\!\!\int_{\Sigma} \overline{\Phi \circ \beta} \, \phi \,\,  n(x)
\cdot {v} \,\, d{v} d\sigma_{\! x} dt
- \int_0^T \!\! \int_\OO \bigl( \Phi \circ \beta ( g ) \, \Lambda_E \, \chi
+  (\Phi \circ \beta)'(g) \, G \, \chi \bigr) \, d{v} dx dt \Bigr| \le \\
&&\qquad\qquad\le [1-\theta + e^{(\theta-1) \, M}] \, C_T,
\eean
and thus
\bean
&&\Bigl| \int_0^T\!\int\!\!\!\int_{\Sigma} \overline{\Phi \circ \beta} \, \phi \,\,  n(x)
\cdot {v} \,\, d{v} d\sigma_{\! x} dt 
- \int_0^T \!\! \int_\OO \bigl[ \Phi \circ \beta ( g ) \, \Lambda_E \, \chi \\
&&\qquad\qquad
+ \bigl( (\Phi \circ \beta)'(g) \, G 
- (\Phi \circ \beta)''(g) \, | \nabla_{v} g|^2 \bigr) \, \chi \bigr] \, d{v} dx dt \Bigr| 
\le 2 \, [1-\theta + e^{(\theta-1) \, M}] \, C_T.
\eean
Once again, by definition of $\gamma g$, we obtain
$$
\Bigl| 
\int_0^T\!\int\!\!\!\int_{\Sigma} \bigl( \overline{\Phi \circ \beta} - 
\Phi \circ \beta ( \gamma \, g ) \bigr) \, \phi \,\,  n(x) \cdot {v} \,\, d{v} d\sigma_{\! x} dt
\Bigr| \le 2 \, [1-\theta + e^{(\theta-1) \, M}] \, C_T \mathop{\longrightarrow}_{M \to \infty,
\theta \to 1} 0,
$$
and $\overline{\Phi \circ \beta} \nearrow \rlim \gamma g_n$ since $\Phi\circ \beta (s) \nearrow s$
when $M \nearrow \infty$, $\theta \searrow 1$, so that $\gamma g = \rlim \gamma g_n$.

\smallskip
In order to prove the a.e. convergence we only have to show, thanks to Proposition~\ref{propWW3}.3, that, up
to the extraction of a subsequence,
\beqn\label{3.21}
\rliminf \beta(\gamma_+ g_n) \ge \beta( \gamma_+ g).
\eeqn
Using Lemma~\ref{Lemma31} and the first step, we can pass to the limit in (\ref{3.17}), up to the extraction of a
subsequence, and we get
\bean
&&\int_0^T\!\!\!\int\!\!\!\int_{\Sigma} \overline{\beta} \, \phi \,\,  n(x)
\cdot {v} \,\, d{v} d\sigma_{\! x} dt
= \int_0^T \!\!\! \int_\OO \bigl(  \beta ( g ) \, \Lambda_{E} \, \chi
+ (\beta'(g) \, G + \beta''(g) \, |\nabla_{v} g|^2 + \mu )\, \chi \bigr) \, d{v} dx dt \\
&&\qquad\qquad= \int_0^T\!\int\!\!\!\int_{\Sigma} (\beta(\gamma g) \,\, n(x) \cdot {v} + \mu ) \, \phi  \,\,
d{v} d\sigma_{\! x} dt ,
\eean
where $\overline{\beta} = \wlim \beta(\gamma g_n)$ is the weak limit in $L^1((0,T) \times
\Sigma)$ of $\beta(\gamma g_n)$. We deduce that $\overline{\beta} \,\, n(x) \cdot {v} = \beta(\gamma
g)
\,\, n(x) \cdot {v} + \mu$ on $(0,T) \times \Sigma$, and in particular
$$
\overline{\beta} \ge \beta(\gamma_+ g)  \quad\hbox{ on } \quad (0,T) \times \Sigma_+.
$$
Since $\rliminf \beta(\gamma_+ g_n)  = \bar\beta$, that ends the proof of (\ref{3.21}).
\qed

%%%%%%%%%%%%%%%%%%%%%%%%%%%%%%%%%%%%%%%%%%%%%%%%
\section{Boltzmann, Vlasov-Poisson and Fokker-Planck equations}\label{sec:Models}
\setcounter{equation}{0}
\setcounter{theo}{0}
%%%%%%%%%%%%%%%%%%%%%%%%%%%%%%%%%%%%%%%%%%%%%%%%

\label{sectionModels}

In this section we derive the a priori physical bound, then make precise the exact meaning of
renormalized solution we deal with and finally state and present a proof of the corresponding stability results. 
In order not to repeat many times the exposition, we consider the full
Vlasov-Poisson-Fokker-Planck-Boltzmann system (VPFPB in short)
\beqn\label{model1}
{\partial f \over \partial t}  +  {v} \cdot \nabla_x  f 
- \hbox{div}_{v} \, \bigl( (\nabla_x V_f + \lambda \,{v}) \, f \bigr) - \nu \, \Delta_{v} f = 
Q(f,f) 
\ \hbox{ in } \ (0,\infty) \times \OO,
\eeqn
where $\nu \ge 0$, $\lambda \in \R$, $Q(f,f)$ stands for the bilinear Boltzmann collision operator  and $V_f$ is given by the mean of the Poisson equation 
\beqn\label{1.7}
- \Delta V_f = \rho_f := \int_{\R^N}Êf \, dv  \ \hbox{ in } \ (0,\infty) \times \Omega, \qquad 
\quad V_f = 0 \ \hbox{ on } \ (0,\infty) \times  \partial\Omega.
\eeqn
We do not give the explicit expression for $Q(f,f)$ that we may find in \cite{26,37} for example. The precise assumptions we make on the cross section are those introduced in \cite{37}. We only recall that the collision operator splits into a gain term and a loss term, $Q(f,f) = Q^+(f,f) - Q^-(f,f)$, that it has the following collision invariants
\beqn\label{model3}
\int_{\R^3} Q(f,f) \, \pmatrix{ 1\cr {v}\cr |{v}|^2\cr } \, d{v} = 0,
\eeqn
and that the so-called {\it entropy production term} $e(f) \ge 0$ satisfies
\beqn\label{model4}
\int_{\R^3} e(f) \, d{v} = - \int_{\R^3} Q(f,f) \, \log f \, d{v}. 
\eeqn
Moreover, it has been established in \cite{37} the following estimate 
\beqn\label{bdQff1}
\forall \, R > 0 \,\,\, \exists C_R < \infty \qquad \int_{B_R} {Q^\pm(f,f) \over 1 + f} \, d{v}
\le C_R  \int_{\R^N} \bigl[ (1+|{v}|^2) \, f + e(f) \bigr] \, d{v},
\eeqn
and in  \cite{63} (we also refer to \cite{50} for a related result) the more accurate estimate
\beqn\label{bdQff2}
\forall \, R > 0 \,\,\, \exists C_R < \infty \qquad \int_{B_R} {|Q(f,f)| \over \sqrt{1 + f}} \, d{v}
\le C_R  \int_{\R^N} \bigl[ (1+|{v}|^2) \, f + e(f) \bigr] \, d{v}. 
\eeqn

We assume furthermore that $f$ satisfies the boundary condition (\ref{BdaryCond}) and the initial condition (\ref{CIf0}), where $f_{in}$ is assumed to verify (\ref{bornef0}), as well as  the following additional bound when $\nu > 0$: 
\beqn\label{nablaV0}
\int_\Omega |\nabla_x V_{f_{in}} |^2 \,  dx < \infty
\quad \hbox{ with } \quad
- \Delta_x V_{f_{in}} = \int_{\R^3} f_{in}(x,{v}) \, d{v} \hbox{ on } \Omega,
\quad V_{f_{in}} = 0 \hbox{ on } \partial\Omega.
\eeqn

\begin{lem}\label{lem:VPFPBound}
For any non-negative initial datum $f_{in}$ such that (\ref{bornef0})-(\ref{nablaV0}) holds
and any time $T \in (0,\infty)$ there exists a constant $C_T \in (0,\infty)$ (only depending on $T$ and on $f_{in}$ through the quantities $C_0$ and $\ |\nabla_x V_{f_{in}} \|_{L^2}$) such that any solution $f$ to the initial boundary value problem (\ref{model1})-(\ref{1.7})-(\ref{BdaryCond})-(\ref{CIf0}) satisfies (at least formally) 
\bear\label{modelBd1}
&&\sup_{[0,T]} \Bigl\{\int \!\!\! \int_\OO f \, \bigl( 1 +  |{v}|^2 + |\log f| \bigr) \, d{v} dx +
\int_\Omega |\nabla_{x} V_f|^2 \, dx  \Bigr\} \\ \nonumber
&& \qquad\qquad\qquad\qquad+ \int_0^T \!\!\! \int \!\!\! \int_\OO \bigl( e(f) + \nu \, {|\nabla_{v} f|^2 \over f} \bigr) \, d{v} dx
dt 
\le C_T, 
\eear
as well as 
\beqn\label{modelBd2}
 \int_0^T \!\int_{\partial\Omega} \Bigl\{ \EE \left( {\gamma_+ f \over M}Ê\right) + \sqrt{\widetilde{\gamap f } \, }  \,\Bigr\} \, d\sigx dt \le C_T,
\eeqn
where $\EE$ is defined in (\ref{defDGI}). It is worth mentioning that the second estimate in (\ref{modelBd2}) is an a posteriori estimate which we deduce from the interior estimate (\ref{modelBd1}) and a Green formula.  
\end{lem}

\noindent{\sl Proof of (\ref{modelBd1}) in Lemma~\ref{lem:VPFPBound}. }
We claim that for $f$ sufficiently regular and decreasing at the infinity  all the integrations (by parts) that we shall perform are allowed. 

\smallskip\noindent
First, we simply integrate the equation (\ref{model1}) over all variables, and we get the conservation of
mass
$$
\int \!\!\! \int_\OO f(t,.)  \, d{v} dx = \int \!\!\! \int_\OO f_{in} \, d{v} dx
\qquad \forall t \ge 0.
$$

\noindent
Next, setting $h_M(s) = s \, \log (s/M)$ and $E=\nabla_x V_f$, we compute
\bean
&& {\partial \over \partial t} h_M(f) +  {v} \cdot \nabla_x h_M(f)
+ \hbox{div}_{v} \bigl( (E + \lambda \, {v}) h_M(f) \bigr) - \nu \, \Delta_{v} h_M(f) 
=  \\
&&\qquad\qquad= h'_M(f) \, Q(f,f) - \nu \, h''_M(f) \, |\nabla_{v} f|^2 - f \, (E + \lambda \, {v}) \cdot
\nabla_{v} (\log M)  \\
&&\qquad\qquad+ \lambda \, (h_M(f) - f \, h'_M(f) ) + 2 \, \nu \, \nabla_{v} f \cdot \nabla_{v} (\log M) + \nu
\, f \, \Delta_{v} (\log M),  
\eean
where $h'_M(s) = 1 + \log \, (s/M)$.
We integrate this equation over the $x,v$ variables using the collision invariants (\ref{model3}) and the entropy production identity (\ref{model4}), to obtain
\bear\label{model5}
&&{d \over dt} \int \!\!\! \int_\OO h_M(f) \, d{v} dx 
+ \int \!\!\! \int_\OO \bigl( e(f) + \nu \, {|\nabla_{v} f|^2 \over f} \bigr) \, d{v} dx
+\int \!\! \int_{\Sigma} h_M(\gamma f ) \,{v} \cdot n(x) \, d{v} d\sigx =\\ \nonumber
&&\qquad\qquad= \int_\Omega  E \cdot { j \over \Theta} \, dx 
+  \int \!\!\! \int_\OO \bigl\{ \lambda \, \bigl( {|{v}|^2  \over \Theta} - 1 ) + 
{\nu \over \Theta} \Bigr\} \, f \, d{v} dx,
\eear
where
$$
j(t,x) = \int_{\R^3} {v} \, f(t,x,{v}) \, d{v}.
$$

\noindent
We first remark that integrating equation (\ref{model1}) in the velocity variable we have
$$
{\partial \over \partial t} \rho + div_x \, j = 0
 \quad \hbox{ on } \ (0,\infty) \times \Omega,
$$
and therefore
\beqn\label{model6}
-\int_\Omega  E \cdot { j \over \Theta} \, dx
= \int_\Omega  \nabla V_f \cdot { j \over \Theta} \, dx  
= \int_\Omega  {V_f  \over \Theta} \, {\partial \rho \over \partial t} \, dx  
= {d \over dt} \int_\Omega  {|\nabla_x V_f |^2 \over 2 \, \Theta} \, dx.
\eeqn
Next, combining (\ref{model5}), (\ref{model6}) and the boundary estimate (\ref{FTbd4}) we obtain
\bean\label{(4.9)}
&&{d \over dt} \Bigl\{\int \!\!\! \int_\OO h_M(f) \, d{v} dx +
\int_\Omega {|\nabla_x V_f |^2  \over 2 \, \Theta} \, dx  \Bigr\} 
+ \int \!\!\! \int_\OO \bigl[ e(f) + \nu \, {|\nabla_{v} f|^2 \over f} \bigr] \, d{v} dx \\
&&\qquad\qquad + \bar\alpha \, \int_{\partial\Omega} \EE ( \gamma_+ \, f) \, d\sigx  
 \le C_{\lambda,\nu} \int \!\!\! \int_\OO (1 + |{v}|^2) \, f \, d{v} dx.
\eean
Here and below, we set $\bar \alpha = \alpha$ in the case of the constant accommodation coefficient (\ref{RLD}) and $\bar\alpha$ is defined just after equation (\ref{NLBC}) in the case of mass flux dependent accommodation coefficient.
Using the elementary estimate (\ref{BddH1}) and (\ref{BddH2}) we conclude that (\ref{modelBd1}) holds, as well as the first estimate in  (\ref{modelBd2}). 

\smallskip
In order to prove the second estimate in  (\ref{modelBd2}), we fix $\chi \in \DD(\R^N)$ such that $0 \le \chi \le 1$, 
$\chi = 1$ on $B_1$ and  supp$\, \chi \subset B_2$ and we apply the Green formula (\ref{3.10})  written with $\phi = n(x)
\cdot {v} \, \chi({v})$ and $\beta(s) = \sqrt{1 + s}$. We get 
\bear\label{sqrtQ1}
&&
\int_0^T \!\!\! \int \!\!\! \int_\Sigma \sqrt{1 +\gamma f}  \, (n(x) \cdot v)^2 \, \chi \, dv  d\sigx dt =
\Bigl[ \, \int\!\!\!\int_{\OO} \sqrt{1 +\gamma f}  \, \phi \,  d{v} dx \, \Bigr]_{T}^{0} \\ \nonumber
&&\qquad
+ \int_{0}^{T}\!\!\int\!\!\!\int_{\OO} \bigl( \sqrt{1 + f} \, \Bigl( v \cdot \nabla_x + (\nabla_x V_f + \lambda \, v) \cdot \nabla_v + \nu \, \Delta_v + N \, \lambda \Bigr) \phi \, dv  dx dt \\ \nonumber
&&\qquad
+ \int_{0}^{T}\!\!\int\!\!\!\int_{\OO} \Bigl( { Q(f,f) \over  2 \, (1 +  f)^{1/2}}  +  { \nu \over  4}  \,  { |\nabla_{v} f|^2\over  (1 +  f)^{3/2}}  \Bigr) \phi \, dv  dx dt .
\eear
Thanks to (\ref{modelBd1}) and (\ref{bdQff2}) and because $\nabla_x \phi \in L^\infty$, $D^2_v \phi \in L^\infty$, we see that the right hand side term in (\ref{sqrtQ1}) is bounded by a constant denoted by $C'_T$ and which only depends on $C_T$ defined in  (\ref{modelBd1}). 
On the other hand, from the boundary condition (\ref{BdaryCond})-(\ref{RLD}) or (\ref{BdaryCond})-(\ref{NLBC}), we have
$\gamam f \ge \bar\alpha \, M({v}) \, \widetilde{\gamap f}$ on $(0,T) \times \Sigma_-$. 
Therefore there is a constant $C_\chi >0$ such that
\bean
C_\chi \int_0^T \!\!\! \int_{\partial\Omega} \sqrt{\widetilde{\gamap f } \, }  \, d\sigx dt
&\le&
\int_0^T \!\!\! \int \!\!\! \int_{\Sigma_-} \sqrt{\widetilde{\gamap f } \, }  \, \bar\alpha^{1/2}
\, M^{1/2}({v}) \, \chi \, (n(x) \cdot v)^2 \, dv  d\sigx dt\\
&\le&
\int_0^T \!\!\! \int \!\!\! \int_{\Sigma_-} \sqrt{\gamam f} \, \chi \, (n(x) \cdot v)^2  \, dv  d\sigx dt\le
C'_T, 
\eean
which ends the proof of (\ref{modelBd2}). 
\qed

\medskip
We can now specify the sense of the solution we deal with.
With DiPerna and Lions \cite{35}, \cite{37,38}, \cite{50} we say that  $0 \le f \in C([0, \infty); L^1(\OO))$ is a
renormalized solution of (\ref{model1})-(\ref{1.7})-(\ref{BdaryCond})-(\ref{CIf0}) if first $f$ satisfies the a priori physical bound
(\ref{modelBd1}) and is a solution of
\bear\label{(4.10)}
&&
{\partial \over \partial t} \beta (f)  + {v} \cdot \nabla_x \beta (f) 
+ (\nabla_x V_f + \lambda \, {v}) \cdot \nabla_{v} \beta (f)  -\nu \, \Delta_{v} \beta(f) = \\ \nonumber
&&\qquad\qquad=\beta '(f) \, (Q(f,f) + \lambda \, N \, f)  - \nu \, \beta''(f) \, |\nabla_{v} f|^2
\quad \hbox{in} \quad {\cal D}'((0,T) \times \OO ), 
\eear
for all time $T>0$, and all $\beta \in \BB_5$, the class of all
functions $\beta \in C^2(\R)$ such that $|\beta''(s)| \le C / (1 + s)$,
$|\beta'(s)| \le C / \sqrt{1 + s}$, $\forall s \ge 0$.  Thanks to (\ref{modelBd1}) (and (\ref{bdQff2})) we see that each term in equation (\ref{(4.10)}) makes sense.
Next, the trace functions $f(0,.)$ and $\gamma f$ defined by Theorem~\ref{theoTraceRenom} through the Green formula
(\ref{3.10})  must satisfy (\ref{CIf0}) and (\ref{BdaryCond}), say almost everywhere. Finally, we will always assume that $\gamma f$ satisfies the additional bound  (\ref{modelBd2}). 

\smallskip\noindent
Our main result is the following stability or compactness result. Once again, in order not to repeat several times the proof, we establish our result for the full VPFPB system and the full VPB system, the same holds for the same equation with less terms. 

\begin{theo}\label{Th:ModStab}
Let $(f_n)$ be a sequence of renormalized solutions to equation (\ref{model1})-(\ref{1.7}) such that the associated trace functions $\gamma f_n$ satisfy (\ref{BdaryCond}), with the linear reflection operator (\ref{RLD}) when $\nu = 0$ and a possibly mass flux depending accommodation coefficient  (\ref{NLBC}) when $\nu > 0$ (FP type models). Let us furthermore assume that both the sequence of solutions $(f_n)$ and  the trace sequence $(\gamma f_n)$ satisfy (uniformly in $n$) the  natural physical a priori bounds 
\bear\label{ModUnifBd1}
&&\sup_{[0,T]} \Bigl\{\int \!\!\! \int_\OO f_n \, \bigl( 1 +  |{v}|^2 + |\log f_n| \bigr) \, d{v} dx +
\int_\Omega |\nabla_{x} V_{f_n}|^2 \, dx  \Bigr\} \\ \nonumber
&& \qquad\qquad+ \int_0^T \!\!\! \int \!\!\! \int_\OO \bigl( e(f_n) + \nu \, {|\nabla_{v} f_n|^2 \over f_n} \bigr) \, d{v} dx
dt + 
 \int_0^T \!\int_{\partial\Omega} \EE \left( {\gamma_+ f_n \over M}Ê\right)  \, d\sigx dt
\le C_T.
\eear
If $f_n(0,.)$ converges to $f_{in}$ weakly in $L^1(\OO)$ then, up to the extraction of a subsequence, $f^n$ converges weakly in $L^p(0,T;L^1(\OO))$ for all $T > 0$ and $p \in [1,\infty)$ (the convergence being strong when $\nu > 0$) to a renormalized solution $f$ to (\ref{model1})-(\ref{1.7})  with initial value $f_{in}$ and which satisfies the physical estimates (\ref{modelBd1}).  
Furthermore, for any $\eps > 0$ and $T > 0$, there exists a measurable set $A \subset (0,T) \times \partial\Omega$ such that 
meas$\, ((0,T) \times \partial\Omega \, \backslash \, A) < \eps$ and 
$$
\gamap f_n \wto\gamap f \quad\hbox{weakly  in}\quad L^1(A \times \R^N,d\lambda_1),
$$
(the convergence being strong when $\nu > 0$). 
As a consequence we can pass to the limit in the boundary reflection condition  (\ref{BdaryCond}) (and (\ref{NLBC}) when $\nu > 0$), so that the trace condition is fulfilled and the trace estimate (\ref{modelBd2}) holds. 
\end{theo}

\noindent
{\sl Proof of the Theorem~\ref{Th:ModStab}.} From (\ref{ModUnifBd1}) we deduce, extracting a subsequence if necessary, that $f_n$ converges weakly in $L^p(0,T;L^1(\OO))$ ($\forall p \in [1,\infty)$) to a function $f$ and that the local mass density $\rho_n = \rho_{f_n}$ satisfies (see \cite{50})
$$
\sup_{[0,T]} \int_\Omega \rho_n (1 + |\log \rho_n|) \, dx \le C_T.
$$

In the case $\nu = 0$, using the velocity averaging lemma of \cite{GLPS,DiLM} and the standard properties of the Poisson equation, we also show  (see for instance \cite{50} and \cite{60}) 
$$
\rho^n \mathop{\longrightarrow}_{n \to \infty} \rho_f \,\hbox{ in } \, L^p(0,T; L^1(\Omega))
\quad\hbox{and}\quad \nabla_x V_{f_n}\mathop{\longrightarrow}_{n \to \infty} \nabla_x V_f
\,\hbox{ in } \, L^p(0,T; W^{1,1} \cap L^a(\Omega))
$$
for all $T \in (0,\infty)$, $p \in [1,\infty)$ and $a \in [1,2)$. It is also shown in \cite{50} that 
$$
{Q^\pm(f_n,f_n) \over 1+\delta \, f^n} \mathop{\rightharpoonup}_{n \to \infty} \bar Q^\pm_\delta
\hbox{ weakly in } L^1((0,T) \times \OO_R) \quad\hbox{and}\quad
 \bar Q^\pm_\delta \nearrow Q^\pm(f,f) \,\, \hbox{a.e.}
$$

In the case $\nu > 0$, since the term on the right hand side of equation (\ref{(4.10)}) is bounded in
$L^1$, thanks to the uniform estimate (\ref{ModUnifBd1}), and since $\Lambda_{E_{f_n}}$ is an hypoelliptic operator (see \cite{35},
\cite{15}, \cite{47}), we obtain that, say,
$\log (1 + f^n)$ and next $f^n$ converge a.e. (see \cite{14} and \cite{35}). We conclude that $f^n \to f$ strongly in
$L^p(0,T;L^1(\OO))$, $\forall p \in [1,\infty)$.  It is also shown in \cite{35} that 
$$
{Q(f_n,f_n) \over 1 + f_n} \to {Q(f,f) \over 1 + f}
\quad\hbox{strongly in}\quad L^1((0,T) \times \OO_R).
$$
Therefore, using Theorem~\ref{th:rcvgeV} or Theorem~\ref{th:rcvgeVFP}, we obtain that $f$ satisfies the
renormalized equation (\ref{(4.10)}) (first for renormalizing function $\beta \in \BB_4$ and next for $\beta \in \BB_5$) and that
$$
\gamma f_n \rto \gamma f \quad\hbox{in the renormalized sense on } (0,T) \times \Sigma,
$$
as well as 
$$
\gamma f_n \rightarrow \gamma f \quad\hbox{a.e. on } (0,T) \times \Sigma,
$$
when $\nu > 0$. 
It is worth mentioning that $f$ also satisfies the physical estimate (\ref{modelBd1}), see \cite{35,38,50}
Next, from (\ref{BdaryCond}) we have
$$
\widetilde{\gamap f_n} \le  \bar\alpha^{-1} \, M^{-1}({v}) \, \gamam f_n 
\quad\hbox{on}\quad (0,T) \times \Sigma_-,
$$
so that 
$$
\widetilde{\gamap f_n} \rto \psi \quad\hbox{in}\quad (0,T) \times \partial\Omega,
\quad\hbox{with}\quad \psi \le \bar\alpha^{-1} \, M^{-1}({v}) \, \gamam f.
$$ 
Furthermore, repeating the proof of Lemma~\ref{lem:VPFPBound} we get that $\psi \in L^{1/2}((0,T) \times
\partial\Omega)$. Now, we can apply Theorem~\ref{theoxi2} (with $m(v) = M(v)$, $y= (t,x)$, $d\varpi_y(v) = {\bf 1}_{\Sigma^x_+} \, |n(x) \cdot v| \, dv$, $\phi_n = \gamma_+ f_n$ and $d\nu(y) = d\sigma_{\! x} dt$),  which says that for every $\eps > 0$ there is
$A=A_\eps \subset (0,T) \times \partial\Omega$ such that meas$\,((0,T) \times \partial\Omega
\backslash A) < \eps$ and 
$$
\gamap f_n \wto \gamap f \quad\hbox{weakly in}\quad L^1( A \times \R^N).
\leqno
$$
In the case $\nu > 0$, since we already know the a.e. convergence, this convergence is in fact strong in $L^1( A \times
\R^N)$. There is no difficulty in passing to the limit in the boundary condition so that $f$
satisfies (\ref{BdaryCond}) and $f$ satisfies the same physical estimate (\ref{modelBd2}) thanks 
to the convexity argument of Theorem \ref{theo1xi}.
\qed

\begin{rem}\label{Neuman} For the Boltzmann equation and the FPB equation, as well as for the VP system
and the VPFP system when the Poisson equation (\ref{1.7}) is provided with Neumann condition, we can
prove the additional a priori estimate (\ref{FTbd2}) on the trace function. 
As a consequence, we may also establish the a priori physical bound (\ref{modelBd1}) for a time and position
dependent wall temperature $\Theta = \Theta(t,x)$ which satisfies
$$
0 < \Theta_0 \le \Theta(t,x) \le \Theta_1 < \infty.
$$
Therefore, the stability result and the corresponding existence result can be generalized to these
kind of boundary conditions. We refer to \cite{6} and \cite{54} for more details.
\end{rem}

\begin{rem}\label{RGal} Consider  the general reflection operator
\beqn\label{(4.16)}
\RR \, \phi = \int_{{v}'Ê\cdot n(x) > 0} k({v},{v}') \, \phi({v}') \, {v}'Ê\cdot n(x) \, d{v}'
\eeqn
where the measurable function $k$ satisfies the usual non-negative, normalization and reciprocity
conditions
\beqn\label{(i)}
k \ge 0, \qquad
\!\!\int_{{v}Ê\cdot n(x) < 0} k({v},{v}') \, d{v} = 1, \qquad
\RR \, M = M, 
\eeqn
where $M$ is the normalized Maxwellian (\ref{defMxi}). For that reflection operator (\ref{(4.16)}), we can prove that a solution $f$ 
to equations (\ref{model1})-(\ref{1.7})-(\ref{BdaryCond}) formally satisfies the a priori physical estimate (\ref{modelBd1})-(\ref{modelBd2}) with $\EE$ replaced by 
$$
\EE_k (\phi/M) := \int_{{v}Ê\cdot n(x) > 0} 
\bigl[ \, h \bigl( {\phi \over M} \bigr) - h \bigl( {\RR \, \phi \over M} \bigr) \, \bigr] 
\, M \, {v}Ê\cdot n(x) \, d{v}.
$$
By Jensen inequality one can prove that $\EE_k$ is non-negative, see \cite{39}, \cite{30}, \cite{41}.
However, we do not know whether our analysis can be adapted to this general kernel. 
Nevertheless, considering a sequence $(f_n)$ of solutions which satisfies the uniform interior estimate in 
(\ref{ModUnifBd1}), we can pass to the limit in (\ref{BdaryCond}) with the help of 
Theorem~\ref{th:rcvgeV} or Theorem~\ref{th:rcvgeVFP} and of Proposition~\ref{propWW3}.4, and we get 
that the limit function $f$ is a solution which trace  $\gamma f$ satisfies the
boundary inequality condition (\ref{bdConRelax}). That extends and generalizes previous results
known for the Boltzmann equation, see for instance \cite{6}, \cite{29}, \cite{54}.
\end{rem}

%%%%%%%%%%%%%%% appendix %%%%%%%%%%%%%%%%%%%%%

\appendix

\section{Appendix: More about the renormalized convergence} \label{sec:A}
\setcounter{equation}{0}
\setcounter{theo}{0}

We come back to the notion of renormalized convergence and mainly discuss its relationship 
 with the biting-$L^1$ weak convergence. 

\begin{rem}\label{rem:Appendix}
1. Hypothesis $\psi \in L^0(Y)$ in Theorem~\ref{theoBit2} (and (\ref{2.8})) is fundamental, since for
example, the sequence $(\psi_n)$ defined by $\psi_n = \psi \equiv +\infty$ $\forall n$ does converge in the
renormalized sense to $\psi$, but $(\psi_n)$ does not converge (and none of its subsequence!) in the biting $L^1$-weak sense. 

2. The (asymptotically) boundedness of $(\psi_n)$ in $L^0$ does not
guarantee that $(\psi_n)$ satisfies, up to the extraction of a subsequence, (\ref{2.7}) or (\ref{2.8}). An instructive example
is the following: we define $u(y) = 1/y$ on $Y = [0,1]$ that we extend by
$1$-periodicity to $\R$, and we set $\psi_n(y) = u(n \, y)$ for $y \in Y$. 
Therefore, $(\psi_n)$ is obviously bounded in $L^a(Y)$ for all $a \in [0,1)$ and converges to
$\psi \equiv +\infty$ in the renormalized sense.
\end{rem}

\begin{prop}\label{prop:Appendix}
1. There exists $(\phi_n)$ which r-converges but does not b-converges.

2. There exists $(\phi_n)$ which b-converges but does not r-converges.

3. Given a sequence $(\phi_n)$, the property 
\beqn\label{6.2}\qquad
\hbox{for any sub-sequence }(\phi_{n'}) \hbox{ there exists a sub-sequence }(\phi_{n''}) \hbox{ of
} (\phi_{n'}) \hbox{ such that } \phi_{n''} \wwto \phi
\eeqn
does not imply $\phi_n \wwto \phi$, where $\wwto$ denotes either the b-convergence or the r-convergence
As a consequence, the b-convergence and the r-convergence are not associated to any Hausdorff (separated) topological structure.
\end{prop}

\noindent
{\sl Proof of Proposition~\ref{prop:Appendix}. } 
{\sl Points 1 \& 3. }ÊLet $(\phi_n)$ be the sequence defined by $\phi_n = \phi_{p,k} = p \, {\bf 1}_{[k/p,(k+1)/p]}$ where $p \in \N^*$, $0 \le k \le p-1$ and  $n = 1 + 2 + ... + p + k$. Then $(\phi_n)$ is bounded in $L^1$ and clearly r-converges to $0$, but does not b-converge.
Moreover, for any subsequence $(\phi_{n'})$ we can find a second subsequence $(\phi_{n''})$ such that $\phi_{n''}$ b-converges to $0$.

\smallskip\noindent
{\sl Points 2 \& 3. } Consider $\mu_y = \mu$ and $\nu_y =\nu$ two Young measures on $Y =
[0,1]$ such that 
\bean
&&\int_\R z \, \mu(dz) = \int_\R z \, \nu(dz) =: \phi \in L^1(Y),\\
&&\int_\R T_{M}(z) \, \mu(dz) \not= \int_\R T_{M}(z) \, \nu(dz) \quad \forall \, M > 0,
\eean
and define $(u_n)$ (resp. $(v_n)$) a sequence of $L^1$ functions tassociated to $\mu$ (resp. $\nu$),
 such that for any $f \in C(\R)$ 
$$
f(u_n) \, \wto \, \bar f := \int_\R f(z) \, \mu_y(dz)
\qquad \bigl( \hbox{resp. } \,\, f(v_n) \, \wto \, \tilde f := \int_\R f(z) \, \nu_y(dz) \bigr), 
$$
see \cite[Theorem 5]{61}, \cite{67}. Then define $(\phi_n)$ by setting $\phi_{2 \, n} = u_n$, $\phi_{2 \, n + 1} = v_n$.
In such a way, we have exhibited  a sequence $(\phi_n)$ which does not r-converge (for instance does not $(T_M)$-renormalized converge) but converges to $\phi$ in the weak $L^1$ sense, and thus b-converges to $\phi$. Moreover, for any sub-sequence 
$(\phi_{n'})$, there exists a second  sub-sequence $(\phi_{n''})$ which either converges to $\bar T_M$ (if $\{Ên' \}$ contain an infinity of even integer numbers) or to $\tilde T_M$ (if $\{Ên' \}$ contain an infinity of odd integer numbers). Because $\bar T_M \nearrow \phi$ and $\tilde T_M \nearrow \phi$ when $M\nearrow \infty$, in both case $\phi_{n''}$ r-converges to $\phi$, and (\ref{6.2}) holds. 
\qed

\bigskip
\noindent
{\bf{Acknowledgments.}}  I would like to thank T. Horsin and O. Kavian for many helpful discussions on weak-weak convergence.
I express my gratitude to J. Soler for his kind hospitality while I was visiting at the Universidad de Granada where part of this work has been done. Financial support from TMR (ERB FMBX CT97 0157) for that visit is gratefully acknowledged.
I also would like to thank F. Murat and F. Castaing for quoting additional references and for their useful comments.

%\newpage

\end{document}